\documentclass[letterpaper, 10 pt, conference]{ieeeconf}%
\usepackage{makeidx}
\usepackage{graphicx}
\usepackage{multicol}
\usepackage[bottom]{footmisc}
\usepackage{latexsym}
\usepackage{subcaption}
\usepackage{cite,color,comment,xspace}
\usepackage{amsfonts}
\usepackage{amsmath}
\usepackage{amssymb}
\usepackage{algorithm}
\usepackage{algorithmic}
\usepackage{url}
\usepackage{times}
\usepackage{mathptmx}
\usepackage{enumerate}
\usepackage{epstopdf}
\usepackage{epsfig}
\usepackage{authblk}
\usepackage[skip=2pt,font=scriptsize]{caption}%
\providecommand{\U}[1]{\protect\rule{.1in}{.1in}}
\renewcommand{\baselinestretch}{.985}
\IEEEoverridecommandlockouts

\newtheorem{theorem}{Theorem}

\makeatletter
\@addtoreset{corollary}{section}

\@addtoreset{lemma}{section}

\@addtoreset{definition}{section}
\makeatother
\allowdisplaybreaks
\begin{document}
\pdfoutput = 1
\title{{\LARGE \textbf{Escaping Local Optima in a Class of Multi-Agent Distributed
Optimization Problems: A Boosting Function Approach}}}
\author{Xinmiao Sun\\Division of Systems Engineering \\Boston University
\and Christos G. Cassandras\\Division of Systems Engineering \\Boston University
\and Kagan Gokbayrak\\Industrial Engineering Department \\Bilkent University\thanks{{\footnotesize The authors' work is supported in
part by NSF under grant CNS-1239021, by AFOSR under grant FA9550-12-1-0113, by
ONR under grant N00014-09-1-1051, and by ARO under Grant W911NF-11-1-0227.
xmsun@bu.edu, cgc@bu.edu, kgokbayr@bilkent.edu.tr}}}
\maketitle

\begin{abstract}
We address the problem of multiple local optima commonly arising in
optimization problems for multi-agent systems, where objective functions are
nonlinear and nonconvex. For the class of coverage control problems, we
propose a systematic approach for escaping a local optimum, rather than
randomly perturbing controllable variables away from it. We show that the
objective function for these problems can be decomposed to facilitate the
evaluation of the local partial derivative of each node in the system and to
provide insights into its structure. This structure is exploited by defining
\textquotedblleft boosting functions\textquotedblright\ applied to the
aforementioned local partial derivative at an equilibrium point where its
value is zero so as to transform it in a way that induces nodes to explore
poorly covered areas of the mission space until a new equilibrium point is
reached. The proposed boosting process ensures that, at its conclusion, the
objective function is no worse than its pre-boosting value. 
However, the global optima cannot be guaranteed. We define three
families of boosting functions with different properties and provide
simulation results illustrating how this approach improves the solutions
obtained for this class of distributed optimization problems.

\end{abstract}

\affil[1]{Division of Systems Engineering, Boston University} \affil[2]{Industrial Engineering Department, Bilkent University}

\thispagestyle{empty} \pagestyle{empty}

\section{Introduction}

Multi-agent systems involve a team of agents (e.g., vehicles, robots, sensor
nodes) that cooperatively perform one or more tasks in a mission space which
may contain uncertainties such as unexpected obstacles or random event
occurrences. The agents communicate, usually wirelessly and over limited
ranges, so there are constraints on the information they can exchange.
Optimization problems are often formulated in the context of such multi-agent
systems and, more often than not, they involve nonlinear, nonconvex objective
functions resulting in solutions where global optimality cannot be easily
guaranteed. The structure of the objective function can sometimes be
exploited, as in cases where it is additive over functions associated with
individual agents; for example, in \cite{MinghuiZhu}, a sum of local nonconvex
objective functions is minimized over nonconvex constraints using an
approximate dual sub-gradient algorithm. In many problems of interest,
however, such an additive structure is not appropriate, as in coverage control
or active sensing \cite{SM2011,CM2004, LM2002,cgc2005} where a set of agents
(typically, sensor nodes) must be positioned so as to cooperatively maximize a
given objective function. In the static version of the problem, the optimal
locations can be determined by an off-line algorithm and nodes will no longer
move. In the dynamic version, nodes may adjust their positions to adapt to
environment changes. Communication costs and constraints imposed on
multi-agent systems, as well as the need to avoid single-point-of-failure
issues, are major motivating factors for developing \emph{distributed}
optimization schemes allowing agents to achieve optimality, each acting
autonomously and with as little information as possible.

\begin{figure}[h]%
\begin{tabular}
[c]{cc}%
\begin{minipage}[t]{1.6in} \includegraphics[
	height=1.15in,
	width=1.35in]
	{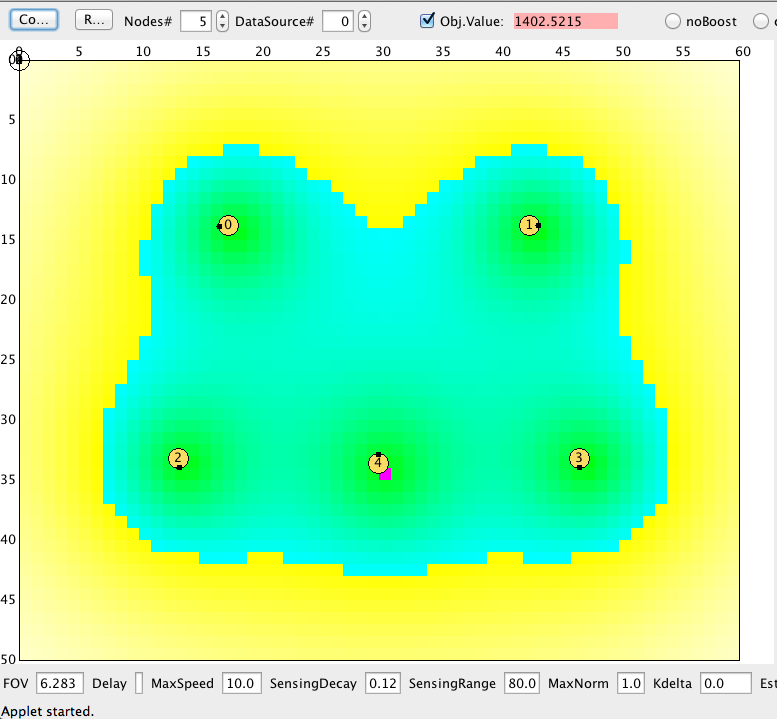} \subcaption{Gradient-based algorithm; \\ optimal obj.function = 1388.1} \label{figObjGradient} \end{minipage} \begin{minipage}[t]{1.6in} \includegraphics[
	height=1.15in,
	width=1.35in]
	{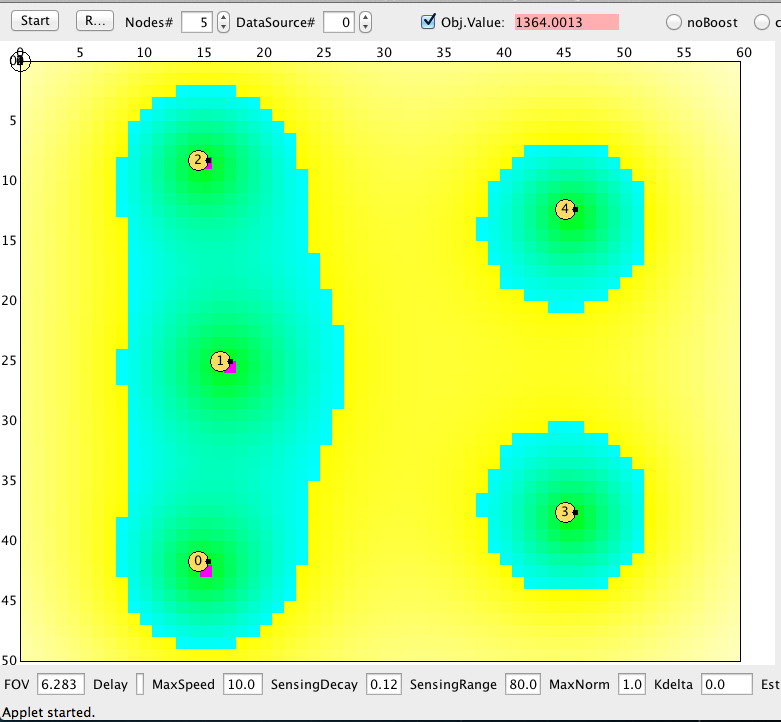} \subcaption{Voronoi patition; \\ optimal obj. function = 1346.5} \label{figObjVoronoi} \end{minipage} 
\end{tabular}
\caption{Comparison between two methods used in a coverage control problem}
\label{figCompare}
\end{figure}

Nonconvex environments for coverage control are treated in
\cite{caicedo2008coverage,caicedo2008performing,breitenmoser2010voronoi,Minyi2011}%
. In \cite{CM2004,gusrialdi2008voronoi,breitenmoser2010voronoi}, algorithms
concentrate on Voronoi partitions of the mission space and the use of Lloyd's
algorithm. We point out that partition-based algorithms do not take into
account the fact that the coverage performance can be improved by sharing
observations made by several nodes. This is illustrated by a simple example in
%Figs. \ref{figObjGradient}-\ref{figObjVoronoi} comparing a common objective
Figure. \ref{figCompare} comparing a common objective
function when a Voronoi partition is used to a distributed gradient-based
approach which optimally positions nodes with overlapping sensor ranges
(darker-colored areas indicate better coverage).

The nonconvexity of objective functions motivates us to seek systematic methods to overcome the
presence of multiple local optima in multi-agent optimization problems. For
off-line centralized solutions, one can resort to global optimization
algorithms that are typically computationally burdensome and time-consuming.
However, for on-line distributed algorithms, this is infeasible; thus, one
normally seeks methods through which controllable variables escape from local
optima and explore the search space of the problem aiming at better
equilibrium points and, ultimately, a globally optimal solution. In
gradient-based algorithms, this is usually done by randomly perturbing
controllable variables away from a local optimum, as in, for example,
simulated annealing \cite{van1987simulated,bertsimas1993simulated} which,
under certain conditions, converges to a global solution in probability.
However, in practice, it is infeasible for agents to perform such a random
search which is notoriously slow and computationally inefficient. \ In the
same vein, in \cite{schwager2008}, a \textquotedblleft ladybug
exploration\textquotedblright\ strategy is applied to an adaptive controller
which aims at balancing coverage and exploration. This approach allows only
two movement directions, thus limiting the ability of agents to explore a
larger fraction of the mission space, especially when obstacles may be
blocking the two exploration directions. In \cite{Minyi2011}, a gradient-based
algorithm was developed to maximize the joint detection probability in a
mission space with obstacles. Recognizing the problem of multiple local
optima, a method was proposed to balance coverage and exploration by modifying
the objective function and assigning a higher reward to points with lower
values of the joint event detection probability metric.

In this paper, we propose a systematic approach for coverage optimization
problems that moves nodes to locations with potentially better performance,
rather than randomly perturbing them away from their current equilibrium. This
is accomplished by exploiting the structure of the problem considered. In
particular, we focus on the class of optimal coverage control problems where
the objective is to maximize the joint detection probability of random events
taking place in a mission space with obstacles. Our first contribution is to
show that each node can decompose the objective function into a local
objective function dependent on this node's controllable position and a
function independent of it. This facilitates the evaluation of the local
partial derivative and provides insights into its structure which we
subsequently exploit. The second contribution is the development of a
systematic method to escape local optima through \textquotedblleft boosting
functions\textquotedblright\ applied to the aforementioned local partial
derivative. The main idea is to alter the local objective function whenever an
equilibrium is reached. A boosting function is a transformation of the
associated local partial derivative which takes place at an equilibrium point,
where its value is zero; the result of the transformation is a non-zero
derivative, which, therefore, forces a node to move in a direction determined
by the boosting function and explore the mission space. When a new equilibrium
point is reached, we revert to the original objective function and the
gradient-based algorithm converges to a new (potentially better and never
worse) equilibrium point. We define three families of boosting functions and
discuss their properties.

In Section II, we formulate the optimization problem and review the
distributed gradient-based solution method developed in \cite{Minyi2011}. In
Section III, we derive the local objective function associated with a node and
its derivative. In Section IV, we introduce the boosting function approach and
three families of boosting functions with different properties. Section V
provides simulation results illustrating how this approach improves the
objective function value and we conclude with Section VI. 

\section{Problem Formulation and Distributed Optimization Solution}

We begin by reviewing the general setting for a large number of multi-agent
control and optimization problems and subsequently concentrate on the optimal
coverage control problem. A \textit{mission space} $\Omega\subset
\mathbb{R}^{2}$ is modeled as a non-self-intersecting polygon, i.e., a polygon
such that any two non-consecutive edges do not intersect. For any $x\in\Omega
$, the function $R(x):\Omega\rightarrow\mathbb{R}$ describes some a priori
information associated with $\Omega$. When the problem is to detect random
events that may take place in $\Omega$, this function captures an a priori
estimate of the frequency of such event occurrences and is referred to as an
\textit{event density} satisfying $R(x)\geq0$ for all $x\in\Omega$ and
$\int_{\Omega}R(x)dx<\infty$. The mission space may contain obstacles modeled
as $m$ non-self-intersecting polygons denoted by $M_{j}$, $j=1,\ldots,m$ which
block the movement of agents. The interior of $M_{j}$ is denoted by
$\mathring{M_{j}}$ and the overall \textit{feasible space} is $F=\Omega
\setminus(\mathring{M_{1}}\cup\ldots\cup\mathring{M_{m}})$, i.e., the space
$\Omega$ excluding all interior points of the obstacles. There are $N$ agents
in the mission space and their positions at time $t$ are defined by $s_{i}%
(t)$, $i=1,\ldots,N$ with an overall position vector $\mathbf{s}%
(t)=(s_{1}(t),\ldots,s_{N}(t))$. Figure. \ref{fig2} shows a mission space with
two obstacles and an agent located at $s_{i}$. The agents may communicate with
each other, but there is generally a limited communication range so that it is
customary to represent such a system as a network of nodes with a link $(i,j)$
defined so that nodes $i,j$ can communicate directly with each other. This
limited communication and the overall cost associated with it are major
motivating factors for developing distributed schemes to allow agents to
operate so as to optimally achieve a given objective with each acting as
autonomously as possible.

In a coverage control problem (e.g., \cite{Minyi2011}%
,\cite{caicedo2008performing},\cite{CM2004}), the agents are \textit{sensor
nodes}. We assume that each such node has a bounded sensing range captured by
the \textit{sensing radius} $\delta_{i}$. Thus, the sensing region of node $i$
is $\Omega_{i}=\{x:d_{i}(x)\leq\delta_{i}\}\text{ where }d_{i}(x)=\Vert x-s_{i}(t)\Vert$.
The presence of obstacles inhibits the sensing ability of a
node, which motivates the definition of a \textit{visibility set}
$V(s_{i})\subset F$ (we omit the explicit dependence of $s_{i}$ on $t$ for
notational simplicity). A point $x\in F$ is \textit{visible} from $s_{i}\in F$
if the line segment defined by $x$ and $s_{i}$ is contained in $F$, i.e.,
$[\lambda x+(1-\lambda)s_{i}]\in F$ for all $\lambda\in\lbrack0,1]$, and $x$
can be sensed, i.e. $x\in\Omega_{i}$. Then, $V(s_{i})=\Omega_{i}\cap\{x:$
$[\lambda x+(1-\lambda)s_{i}]\in F\}$ is a set of points in $F$ which are
visible from $s_{i}$. We also define $\bar{V}(s_{i})=F\setminus V(s_{i})$
to be the \emph{invisibility set} (e.g., the grey area in Fig. \ref{fig2}).

A sensing model for any node $i$ is given by the probability that $i$ detects
an event occurring at $x\in V(s_{i})$, denoted by $p_{i}(x,s_{i})$. We assume
that $p_{i}(x,s_{i})$ is expressed as a function of $d_{i}(x)=\Vert
x-s_{i}\Vert$ and is monotonically decreasing and differentiable in $d_{i}%
(x)$. An example of such a function is $p_{i}(x,s_{i})=p_{0i}e^{-\lambda
_{i}\Vert x-s_{i}\Vert}$. For points that are invisible by node $i$, the
detection probability is zero. Thus, the overall \textit{sensing detection
probability} is denoted as $\hat{p}_{i}(x,s_{i})$ and defined as
\begin{equation}
{\label{SensingModel}}\hat{p}_{i}(x,s_{i})=%
\begin{cases}
p_{i}(x,s_{i}) & \text{if}\quad x\in V(s_{i})\\
0 & \text{if}\quad x\in\bar{V}(s_{i})
\end{cases}
\end{equation}
Note that $\hat{p}_{i}(x,s_{i})$ is not a continuous function of $s_{i}$. We
may now define the \textit{joint detection probability} that an event at
$x\in\Omega$ is detected by at least one of the $N$ cooperating nodes in the
network:
\begin{equation}
{\label{jointP}}P(x,\mathbf{s})=1-\prod_{i=1}^{N}[1-\hat{p}_{i}(x,s_{i})]
\end{equation}
where we have assumed that detection events by nodes are independent. Finally,
assuming that $R(x)=0$ for $x\notin F$, we define the optimal coverage control
problem to maximize $H(\mathbf{s})$, where
\begin{equation}
{\label{obj}}%
\begin{split}
H(\mathbf{s}) &  =\int_{F}R(x)P(x,\mathbf{s})dx\\
s.t.\quad s_{i} &  \in F,i=1,\ldots,N
\end{split}
\end{equation}
Thus, we seek to control the node position vector $\mathbf{s}=(s_{1}%
,\ldots,s_{N})$ so as to maximize the overall joint detection probability of
events taking place in the environment. Note that this is a nonlinear,
generally nonconvex, optimization problem and the objective function
$H(\mathbf{s})$ cannot be expressed in an additive form such as $\sum
_{i=1}^{N}H_{i}(\mathbf{s})$.

As already mentioned, it is highly desirable to develop distributed
optimization algorithms to solve (\ref{obj}) so as to $(i)$ limit costly
communication among nodes (especially in wireless settings where it is known
that communication consumes most of the energy available at each node relative
to sensing or computation tasks) and $(ii)$ impart robustness to the system as
a whole by avoiding single-point-of-failure issues. Towards this goal, a
distributed gradient-based algorithm was developed in \cite{Minyi2011} based
on the iterative scheme:
\begin{equation}
{\label{mov}}s_{i}^{k+1}=s_{i}^{k}+\zeta_{k}\frac{\partial H(\mathbf{s}%
)}{\partial s_{i}^{k}},\text{ \ }k=0,1,\ldots
\end{equation}
where the step size sequence $\{{\zeta_{k}\}}$ is appropriately selected (see
\cite{Bertsekas95}) to ensure convergence of the resulting node trajectories.
If nodes are mobile, then (\ref{mov}) can be interpreted as a motion control
scheme for the $i$th node. In general, a solution through (\ref{mov}) can only
lead to a local maximum and it is easy to observe that many such local maxima
result in poor performance \cite{Minyi2011} (we will show such examples in
Section V).

Our approach in what follows is to first show that $H(\mathbf{s})$ can be
decomposed into a \textquotedblleft local objective function\textquotedblright%
\ $H_{i}(\mathbf{s})$ and a function independent of $s_{i}$ so that node $i$
can locally evaluate its partial derivative with respect to its own
controllable position through $H_{i}(\mathbf{s})$ alone. Our idea then is to
alter $H_{i}(\mathbf{s})$ after a local optimum is attained when
$\frac{\partial H_{i}(\mathbf{s})}{\partial s_{i}}=0$, and to define a new
objective function $\hat{H}_{i}(\mathbf{s})$. By doing so, we force
$\frac{\partial\hat{H}_{i}(\mathbf{s})}{\partial s_{i}}\neq0$, therefore, node
$i$ can \textquotedblleft escape\textquotedblright\ the local optimum and
explore the rest of the mission space in search of a potentially better
equilibrium point. Because of the structure of $\frac{\partial H_{i}%
(\mathbf{s})}{\partial s_{i}}$ and the insights it provides, however, rather
than explicitly altering $H_{i}(\mathbf{s})$ we instead alter $\frac{\partial
H_{i}(\mathbf{s})}{\partial s_{i}}$ through what we refer to in Section IV as
a \emph{boosting function}.

\section{Local Objective Functions for Distributed Gradient-based Algorithms}

We begin by defining $B_{i}$ to be a set of 
 nodes with respect to
$i$:
\begin{equation}
B_{i}=\{k:\Vert s_{i}-s_{k}\Vert<2\delta_{i},\text{ }k=1,\ldots N,\text{
}k\neq i\}\label{Bi}%
\end{equation}
Clearly, this set includes all nodes $k$ whose sensing region $\Omega_{k}$ has
a nonempty intersection with $\Omega_{i}$, the sensing region of node $i$.
Accordingly, given that there is a total number of $N$ nodes, we define a
complementary set $C_{i}$
\begin{equation}
C_{i}=\{k:k\notin B_{i},\text{ }k=1,\ldots N,\text{ }k\neq i\}
\end{equation}
In addition, let $\Phi_{i}(x)$ denote the joint probability that a point
$x\in\Omega$ is not detected by any neighbor node of $i$, defined as
\begin{equation}
{\label{phi}}\Phi_{i}(x)=\prod_{k\in B_{i}}[1-\hat{p}_{k}(x,s_{k})]
\end{equation}
Similarly, let $\bar{\Phi}_{i}(x)$ denote the probability that a point
$x\in\Omega$ is not covered by nodes in $C_{i}$:
\begin{equation}
\bar{\Phi}_{i}(x)=\prod_{j\in C_{i}}[1-\hat{p}_{j}(x,s_{j})]
\end{equation}
The following theorem establishes the decomposition of $H(\mathbf{s})$ into a
function dependent on $s_{i}$, for any $i=1,\ldots,N$, and one dependent on
all other node positions except $s_{i}$.

\begin{theorem}
The objective function $H(\mathbf{s})$ can be written as:%
\begin{equation}
H(\mathbf{s})=H_{i}(\mathbf{s})+\tilde{H}(\bar{\mathbf{s}}_{i})\label{thm}%
\end{equation} 
for any $i=1, \ldots, N$, 
where $\bar{\mathbf{s}}_{i}=(s_{1},\ldots,s_{i-1},s_{i+1},\ldots,s_{N})$, and
\begin{align*}
H_{i}(\mathbf{s}) &  =\int_{V(s_{i})}R(x)\Phi_{i}(x)p_{i}(x,s_{i})dx\\
\tilde{H}(\bar{\mathbf{s}}_{i}) & = \int_{F}R(x)\lbrace 1-\prod_{k=1, k\neq i}^{N}[1-\hat{p}_k(x,s_k)]\rbrace dx
\end{align*}
\end{theorem}

%\begin{proof}
\emph{Proof: }Since $F=V(s_{i})\cup\bar{V}(s_{i})$ and $V(s_{i})\cap\bar
{V}(s_{i})=\varnothing$, we can rewrite $H(\mathbf{s})$ in (\ref{obj}) as the
sum of two integrals:
\begin{equation}
{\label{separateH}}H(\mathbf{s})=\int_{V(s_{i})}R(x)P(x,\mathbf{s}%
)dx+\int_{\bar{V}(s_{i})}R(x)P(x,\mathbf{s})dx
\end{equation}
which we will refer to as $H_{i}^{1}(\mathbf{s})$ and $H_{i}^{2}(\mathbf{s})$,
respectively. Using the definitions of $\Phi_{i}(x)$ and $\bar{\Phi}_{i}(x)$,
the joint detection probability $P(x,\mathbf{s})$ in (\ref{jointP}) can be
written as
\begin{equation}
{\label{P}}P(x,\mathbf{s})=1-\Phi_{i}(x)\bar{\Phi}_{i}(x)[1-p_{i}(x)]
\end{equation}
The integral domain in $H_{i}^{1}(\mathbf{s})$ is the visible set for $s_{i}$
and, from (\ref{SensingModel}) we have $p_{i}(x,s_{i})\neq0$ and
$p_{j}(x,s_{j})=0$ for $j\in C_{i}$, hence, $\bar{\Phi}_{i}(x)=1$. Thus,
$H_{i}^{1}(\mathbf{s})$ can be written as
\begin{equation}
{\label{Hi1}}%
\begin{split}
&  H_{i}^{1}(\mathbf{s})=\int_{V(s_{i})}R(x)[1-\Phi_{i}(x)(1-p_{i}(x))]dx\\
&  =\int_{V(s_{i})}R(x)p_{i}(x,s_{i})\Phi_{i}(x)dx+\int_{V(s_{i})}%
R(x)[1-\Phi_{i}(x)]dx\\
&  =\int_{V(s_{i})}R(x)p_{i}(x,s_{i})\Phi_{i}(x)dx+\int_{V(s_{i})}%
R(x)[1-\Phi_{i}(x)\bar{\Phi}_{i}(x)]dx
\end{split}
\end{equation}
For the $H_{i}^{2}(\mathbf{s})$ term, the integral domain is the invisible set
of $s_{i}$, which implies that $p_{i}(x,s_{i})=0$ for $x\in\bar{V}(s_{i})$.
Using the form of $P(x,\mathbf{s})$ defined in (\ref{P}), $H_{i}%
^{2}(\mathbf{s})$ can be written as%
\begin{equation}
{\label{Hi2}}H_{i}^{2}(\mathbf{s})=\int_{\bar{V}(s_{i})}R(x)[1-\Phi_{i}%
(x)\bar{\Phi}_{i}(x)]dx
\end{equation}
Combining (\ref{Hi1}) and (\ref{Hi2}) and merging the second integral in
(\ref{Hi1}) with the integral in (\ref{Hi2}), we obtain:
\begin{equation*}
\begin{split}
H(\mathbf{s})&=\int_{V(s_{i})}R(x)\Phi_{i}(x)p_{i}(x,s_{i})dx+\int_{F}R(x)[1-\Phi
_{i}(x)\bar{\Phi}_{i}(x)]dx \\
&=\int_{V(s_{i})}R(x)\Phi_{i}(x)p_{i}(x,s_{i})dx \\
&+\int_{F}R(x)\big [1-\prod_{k=1,k\neq i}^{N}[1-\hat{p}_k(x,s_k)]\big ] dx
\end{split}
\end{equation*}
The first term is dependent on $s_{i}$, while the second term is independent
of $s_{i}$ in both integrand and integral domain. Using $\bar{\mathbf{s}}%
_{i}=(s_{1},\ldots,s_{i-1},s_{i+1},\ldots,s_{N})$ to denote a vector of all
node positions except $i$, we define $H_{i}(\mathbf{s})$ and $\tilde{H}%
(\bar{\mathbf{s}}_{i})$ as%
\begin{align*}
H_{i}(\mathbf{s}) &  =\int_{V(s_{i})}R(x)\Phi_{i}(x)p_{i}(x,s_{i})dx\\
\tilde{H}(\bar{\mathbf{s}}_{i}) & = \int_{F}R(x)\lbrace 1-\prod_{k=1, k\neq i}^{N}[1-\hat{p}_k(x,s_k)]\rbrace dx
\end{align*}
and the result follows. $\blacksquare$
%\end{proof}
%=\int_{F}R(x)[1-\Phi_{i}(x)\bar{\Phi}_{i}(x)]dx \\

We refer to $H_{i}(\mathbf{s})$ as the \emph{local objective function} of node
$i$ and observe that it depends on $V(s_{i})$, $p_{i}(x,s_{i})$, and $\Phi
_{i}(x)$ which are all available to node $i$ (the latter through some
communication with nodes in $B_{i}$). This result enables a distributed 
gradient-based optimization solution approach
with each node evaluating $\frac{\partial H_{i}(\mathbf{s})}{\partial s_{i}}$.
We now proceed to derive this derivative using the same method as in
\cite{cassandras2008distributed}. Based on the extension of the Leibnitz rule
\cite{flanders}, we get
\begin{equation}
{\label{derivativeH}}%
\begin{split}
\frac{\partial H_{i}(\mathbf{s})}{\partial s_{ix}} &  =\frac{\partial
}{\partial s_{ix}}\int_{V(s_{i})}R(x)\Phi_{i}(x)p_{i}(x,s_{i})dx\\
&  =\int_{V(s_{i})}R(x)\Phi_{i}(x)\frac{\partial p_{i}(x,s_{i})}{\partial
s_{ix}}dx\\
&  +\int_{\partial V(s_{i})}R(x)\Phi_{i}(x)p_{i}(x,s_{i})(u_{x}dx_{y}%
-u_{y}dx_{x})
\end{split}
\end{equation}
where ($u_{x},u_{y}$) illustrates the \textquotedblleft
velocity\textquotedblright\ vector at a boundary point $x=(x_{x},x_{y})$ of
$V(s_{i})$. The first term, denoted by $E_{ix}$, is
\begin{equation}
{\label{derivativeHfirstpartx}}%
\begin{split}
E_{ix} &  =\int_{V(s_{i})}R(x)\Phi_{i}(x)\frac{\partial p_{i}(x,s_{i}%
)}{\partial s_{ix}}dx\\
&  =\int_{V(s_{i})}R(x)\Phi_{i}(x)\left[  -\frac{dp_{i}(x,s_{i})}{dd_{i}%
(x)}\right]  \frac{(x-s_{i})_{x}}{d_{i}(x)}dx
\end{split}
\end{equation}
where $(x-s_{i})_{x}$ is the $x$ component of the vector $(x-s_{i})$.
Similarly, we can obtain an integral $E_{iy}$ with $(x-s_{i})_{y}$ in place of
$(x-s_{i})_{x}$.
%\begin{equation}{\label{derivativeHfirstparty}}
%\begin{split}
%E_{iy}&=\int_{V(s_i)} R(x)\Phi_i(x)\frac{\partial p_i(x,s_i)}{\partial s_{ix}}dx \\
%&=\int_{V(s_i)}R(x)\prod_{k \in B_i} [1-\hat{p}_k(x,s_k)]\frac{dp_i(x,s_i)}{dd_i(x)}\frac{(s_i-x)_y}{d_i(x)}dx
%\end{split}
%\end{equation}

Let $E_{i}=(E_{ix},E_{iy})$. The integrand of $E_{i}$ can be viewed as a
weighted normalized direction vector $\frac{(x-s_{i})}{d_{i}(x)}$ connecting
$s_{i}$ to $x\in F$ where $x$ is visible by the $i$th node. This weight is
defined as
\begin{equation}
{\label{weight1}}w_{1}(x,\mathbf{s})=-R(x)\Phi_{i}(x)\frac{dp_{i}(x,s_{i})}%
{dd_{i}(x)}%
\end{equation}
Observe that $w_{1}(x,\mathbf{s})\geq0$ because $\frac{dp_{i}(x,s_{i})}{dd_{i}%
(x)}<0$ since $p_{i}(x,s_{i})$ is a decreasing function of $d_{i}$.

\begin{figure}[ptb]
\centering
\includegraphics[
height=2.0in,
width=2.3in]
{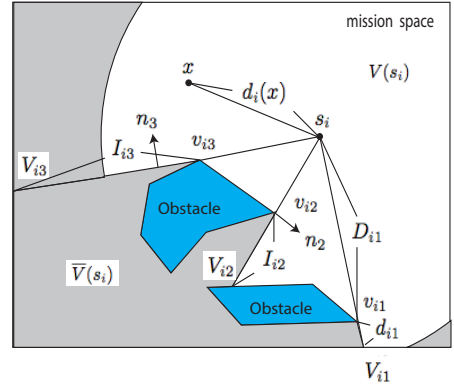} \caption{Mission space with two polygonal obstacles}%
\label{fig2}%
\end{figure}

Next, we evaluate the second term in (\ref{derivativeH}), referred to as
$E_{b}$. This evaluation is more elaborate and requires some additional
notation (see Fig. \ref{fig2}). 
Let $v$ be a reflex vertex(definition can be found in \cite{Minyi2011}) of an obstacle and
let $x\in F$ be a point visible from $v$. A set of points $I(v,x)$, which is a
ray starting from $v$ and extending in the direction of $v-x$, is defined by
\begin{equation}
{\label{line}}I(v,x)=\{q\in V(v):q=\lambda v+(1-\lambda)x,\lambda>1\}
\end{equation}
The ray intersects the boundary of $F$ at an \textit{impact} point. The line
from $v$ to the impact point is a $I(v,x)$.

An \textit{anchor} of $s_{i}$ is a reflex vertex $v$ such that it is visible
from $s_{i}$ and $I(v,s_{i})$ defined in (\ref{line}) is not empty. Denote the
anchors of $s_{i}$ by $v_{ij}$, $j=1,\dots,Q(s_{i})$, where $Q(s_{i})$ is the
number of anchors of $s_{i}$. An \textit{impact} point of $v_{ij}$, denoted by
$V_{ij}$, is the intersection of $I(v_{ij},s_{i})$ and $\partial F$. As an
example, in Fig. \ref{fig2}, $v_{i1}$, $v_{i2}$, $v_{i3}$ are anchors of
$s_{i}$, and $V_{i1}$, $V_{i2}$, $V_{i3}$ are the corresponding impact points.
Let $D_{ij}=\Vert s_{i}-v_{ij}\Vert$ and $d_{ij}=\Vert V_{ij}-v_{ij}\Vert$.
Define $\theta_{ij}$ to be the angle formed by $s_{i}-v_{ij}$ and the
$x$-axis, which satisfies $\theta_{ij}\in\lbrack0,\pi/2]$, that is,
$
\theta_{ij}=arctan \frac{|s_{i}-v_{ij}|_{y}}{|s_{i}-v_{ij}|_{x}}
$.
Using this notation, a detailed derivation of the second term in
(\ref{derivativeH}) may be found in \cite{cassandras2008distributed} with the
final result being:
\begin{equation}
{\label{lineintresultx}}%
\begin{split}
&  E_{bx}=\\
&  \sum_{j\in\Gamma_{i}}sgn(n_{jx})\frac{sin\theta_{ij}}{D_{ij}}\int
_{0}^{z_{ij}}R(\rho_{ij}(r))\Phi_{i}(\rho_{ij}(r))p_{i}(\rho_{ij}(r),s_{i})rdr
\end{split}
\end{equation}
where $\Gamma_{i}=\{j:D_{ij}<\delta_{i},j=1,\ldots,Q(s_{i})\}$; $z_{ij}%
=\min(d_{ij},\delta_{i}-D_{ij})$ and $\rho_{ij}(r)$ is the Cartesian
coordinate of a point on $I_{ij}$ which is a distance $r$ from $v_{ij}$:
\begin{equation}
\rho_{ij}(r)=(V_{ij}-v_{ij})\frac{r}{d_{ij}}+v_{ij}%
\end{equation}
In the same way, we can also obtain $E_{by}$. Note that $E_{b}=(E_{bx}%
,E_{by})$ is the gradient component in (\ref{derivativeH}) due to points on
the boundary $\partial V(s_{i})$. In particular, for each boundary, this
component attracts node $i$ to move in a direction perpendicular to the
boundary and pointing towards $V(s_{i})$. $\ $We can see in
(\ref{lineintresultx}) that every point $x$ written as $\rho_{ij}(r)$ in the
integrand has an associated weight which we define as $w_{2}(x,\mathbf{s})$:
\begin{equation}
{\label{weight2}}w_{2}(x,\mathbf{s})=R(x)\Phi_{i}(x)p_{i}(x,s_{i})
\end{equation}
and observe that $w_{2}(x,\mathbf{s})\geq0$, as was the case for $w_{1}(x,\mathbf{s})$.
Combining (\ref{derivativeHfirstpartx}) and (\ref{lineintresultx}) we finally
obtain the derivative of $H_{i}(\mathbf{s})$ with respect to $s_{i}$:
\begin{equation}
{\label{parx}}%
\begin{split}
&  \frac{\partial{H_{i}(\mathbf{s})}}{\partial{s_{ix}}}=\int_{V(s_{i}%
)}R(x)\prod_{k\in B_{i}}[1-\hat{p}_{k}(x,s_{k})]\frac{dp_{i}(x,s_{i})}%
{dd_{i}(x)}\frac{(s_{i}-x)_{x}}{d_{i}(x)}dx+\\
&  \sum_{j\in\Gamma_{i}}sgn(n_{jx})\frac{sin\theta_{ij}}{D_{ij}}\int
_{0}^{z_{ij}}R(\rho_{ij}(r))\Phi_{i}(\rho_{ij})(r)p_{i}(\rho_{ij}(r),s_{i})rdr
\end{split}
\end{equation}%
\begin{equation}
{\label{pary}}%
\begin{split}
&  \frac{\partial{H_{i}(\mathbf{s})}}{\partial{s_{iy}}}=\int_{V(s_{i}%
)}R(x)\prod_{k\in B_{i}}[1-\hat{p}_{k}(x,s_{k})]\frac{dp_{i}(x,s_{i})}%
{dd_{i}(x)}\frac{(s_{i}-x)_{y}}{d_{i}(x)}dx+\\
&  \sum_{j\in\Gamma_{i}}sgn(n_{jy})\frac{cos\theta_{ij}}{D_{ij}}\int
_{0}^{z_{ij}}R(\rho_{ij}(r))\Phi_{i}(\rho_{ij})(r)p_{i}(\rho_{ij}(r),s_{i})rdr
\end{split}
\end{equation}
We observe that $\frac{\partial{H_{i}(\mathbf{s})}}{\partial{s_{ix}}}$,
$\frac{\partial{H_{i}(\mathbf{s})}}{\partial{s_{iy}}}$ in (\ref{parx}%
)-(\ref{pary}) are the same as $\frac{\partial{H(\mathbf{s})}}{\partial
{s_{ix}}}$, $\frac{\partial{H(\mathbf{s})}}{\partial{s_{iy}}}$, the partial
derivatives of the original objective function ${H(\mathbf{s})}$ which was
derived in \cite{cassandras2008distributed}. In other words, $\frac
{\partial{H(\mathbf{s})}}{\partial{s_{i}}}=\frac{\partial{H_{i}(\mathbf{s})}%
}{\partial{s_{i}}}$, confirming (as expected) that the local objective
function $H_{i}(\mathbf{s})$ is sufficient to provide the required derivative
for a distributed gradient-based algorithm using (\ref{mov}). As pointed out
in \cite{cassandras2008distributed}, the derivation of (\ref{parx}%
)-(\ref{pary}) excludes pathological cases where $s_{i}$ coincides with a
reflex vertex, a polygonal inflection, or a bitangent, where ${H(\mathbf{s})}$
is generally not differentiable.

We can now use the weight definitions (\ref{weight1}) and (\ref{weight2}) in
(\ref{parx}) and (\ref{pary}) to obtain the following form of the local
derivative evaluated by node $i$:
\begin{equation}
{\label{weightedparx}}%
\begin{split}
\frac{\partial{H_{i}(\mathbf{s})}}{\partial{s_{ix}}}= &  \int_{V(s_{i})}%
w_{1}(x,\mathbf{s})\frac{(x-s_{i})_{x}}{d_{i}(x)}dx\\
+ &  \sum_{j\in\Gamma_{i}}sgn(n_{jx})\frac{sin\theta_{ij}}{D_{ij}}\int
_{0}^{z_{ij}}w_{2}(\rho_{ij}(r),s_{i})rdr
\end{split}
\end{equation}%
\begin{equation}
{\label{weightedpary}}%
\begin{split}
\frac{\partial{H_{i}(\mathbf{s})}}{\partial{s_{iy}}}= &  \int_{V(s_{i})}%
w_{1}(x,\mathbf{s})\frac{(x-s_{i})_{y}}{d_{i}(x)}dx\\
+ &  \sum_{j\in\Gamma_{i}}sgn(n_{jy})\frac{cos\theta_{ij}}{D_{ij}}\int
_{0}^{z_{ij}}w_{2}(\rho_{ij}(r),s_{i})rdr
\end{split}
\end{equation}
We can see that the essence of each derivative is captured in the weights
$w_{1}(x,\mathbf{s}),$ $w_{2}(x,\mathbf{s})$. In the first integral, $w_{1}(x,\mathbf{s})$
controls the mechanism through which node $i$ is attracted to different points
$x\in V(s_{i})$ through $\frac{(x-s_{i})}{d_{i}(x)}$. If obstacles are
present, then $w_{2}(x,\mathbf{s})$ in the second integral controls the attraction
that boundary points exert on node $i$ with the geometrical features of the
mission space contributing through $n_{jx}$, $n_{jy}$, $\theta_{ij}$, and
$D_{ij}$ in (\ref{weightedparx})-(\ref{weightedpary}). This viewpoint
motivates the boosting function approach described next.

\section{The Boosting Function Approach}

As defined in (\ref{obj}), this nonlinear, generally nonconvex, optimization
problem may have multiple local optima to which a gradient-based algorithm may
converge. When we apply a distributed optimization algorithm based on
$\frac{\partial{H_{i}(\mathbf{s})}}{\partial{s_{i}}}$ as described above, any
equilibrium point is characterized by $\frac{\partial{H_{i}(\mathbf{s})}%
}{\partial{s_{i}}}=0$. Since node $i$ controls its position based on its local
objective function ${H_{i}(\mathbf{s})}$, a simple way to \textquotedblleft
escape\textquotedblright\ a local optimum ${\mathbf{s}}^{1}$ is to alter
${H_{i}(\mathbf{s})}$ by replacing it with some ${\hat{H}_{i}(\mathbf{s})\neq
H_{i}(\mathbf{s})}$ thus forcing $\left.  \frac{\partial{\hat{H}%
_{i}(\mathbf{s})}}{\partial{s_{i}}}\right\vert _{s_{i}^{1}}\neq0$ and inducing
the node to explore the rest of the mission space for potentially better
equilibria. Subsequently, when a new equilibrium is reached with node $i$ at
$\tilde{s}_{i}^{1}\neq s_{i}^{1}$ and $\left.  \frac{\partial{\hat{H}%
_{i}(\mathbf{s})}}{\partial{s_{i}}}\right\vert _{\tilde{s}_{i}^{1}}=0$, we can
revert to ${H_{i}(\mathbf{s})}$, which, in turn will force $\left.
\frac{\partial{H_{i}(\mathbf{s})}}{\partial{s_{i}}}\right\vert _{\tilde{s}%
_{i}^{1}}\neq0$ and the node will seek a new equilibrium at $s_{i}^{2}$.

Selecting the proper ${\hat{H}_{i}(\mathbf{s})}$ to temporarily replace
${H_{i}(\mathbf{s})}$ is not a simple process. However, focusing on
$\frac{\partial{H_{i}(\mathbf{s})}}{\partial{s_{i}}}$ instead of
${H_{i}(\mathbf{s})}$ is much simpler due to the nature of the derivatives we
derived in (\ref{weightedparx})-(\ref{weightedpary}). In particular, the
effect of altering ${H_{i}(\mathbf{s})}$ can be accomplished by transforming
the weights $w_{1}(x,\mathbf{s})$, $w_{2}(x,\mathbf{s})$ in (\ref{weightedparx}%
)-(\ref{weightedpary}) by \textquotedblleft boosting\textquotedblright\ them
in a way that forces $\frac{\partial{H_{i}(\mathbf{s})}}{\partial{s_{i}}}=0$
at a local optimum to become nonzero. The net effect is that the attraction
exerted by some points $x\in F$ on $s_{i}$ is \textquotedblleft
boosted\textquotedblright\ so as to promote exploration of the mission space
by node $i$ in search of better optima.

In contrast to various techniques which aim at randomly perturbing
controllable variables away from a local optimum (e.g., simulated annealing),
this approach provides a systematic mechanism for accomplishing this goal by
exploiting the structure of the specific optimization problem reflected
through the form of the derivatives (\ref{weightedparx})-(\ref{weightedpary}).
Specifically, it is clear from these expressions that this can be done by
assigning a higher weight (i.e., boosting) to directions in the mission space
that provide greater opportunity for exploration and, ultimately
\textquotedblleft better coverage\textquotedblright. To develop such a
systematic approach, we define transformations of the weights $w_{1}(x,s_{i}%
)$, $w_{2}(x,\mathbf{s})$ for interior points and for boundary points respectively
as follows:
\begin{align}
\hat{w}_{1}(x,\mathbf{s}) & =g_{i}(w_{1}(x,\mathbf{s})) \\
\hat{w}_{2}(x,\mathbf{s}) & =h_{i}(w_{2}(x,\mathbf{s})) 
\end{align}
where $g_{i}(\cdot)$ and $h_{i}(\cdot)$ are functions of the original weights
$w_{1}(x,\mathbf{s})$ and $w_{2}(x,\mathbf{s})$ respectively. We refer to $g_{i}(\cdot)$
and $h_{i}(\cdot)$ as \emph{boosting functions} for node $i=1,\ldots,N$. Note
that these may be node-dependent and that each node may select the time at
which this boosting is done, independent from other nodes. In other words, the
boosting operation may also be implemented in distributed fashion, in which
case we refer to this process at node $i$ as \emph{self-boosting}.

In the remainder of this paper, we concentrate on functions $g_{i}(\cdot)$ and
$h_{i}(\cdot)$ which have the form
\begin{align}
{\label{modifiedlineari}}\hat{w}_{1}(x,\mathbf{s}) &=\alpha_{i1}(x,\mathbf{s}%
)w_{1}(x,\mathbf{s})+\beta_{i1}(x,\mathbf{s}) \\
{\label{modifiedlinearb}}\hat{w}_{2}(x,\mathbf{s}) &=\alpha_{i2}(x,\mathbf{s}%
)w_{2}(x,\mathbf{s})+\beta_{i2}(x,\mathbf{s})
\end{align}
where $\alpha_{i1}(x,\mathbf{s})$, $\beta_{i1}(x,\mathbf{s})$, $\alpha
_{i2}(x,\mathbf{s})$, and $\beta_{i2}(x,\mathbf{s})$ are functions dependent
on the point $x$ and the node position vector $\mathbf{s}$ in general. We
point out that although the form of (\ref{modifiedlineari}%
)-(\ref{modifiedlinearb}) is linear, the functions $\alpha_{ij}(x,\mathbf{s}%
)$, $\beta_{ij}(x,\mathbf{s})$, $j=1,2$, $i=1,\ldots,N$ are generally
nonlinear in their arguments.

To keep notation simple, let us concentrate on a single node $i$ and omit the
subscript $i$ in $\alpha_{ij}(x,\mathbf{s})$, $\beta_{ij}(x,\mathbf{s})$
above. By replacing $w_{1}(x,\mathbf{s})$, $w_{2}(x,\mathbf{s})$ with $\hat{w}%
_{1}(x,s_{i})$, $\hat{w}_{2}(x,s_{i})$ respectively, we obtain the
\emph{boosted derivative} $\frac{\partial{\hat{H}(\mathbf{s})}}{\partial
s_{i}}$ as follows 
%(we limit ourselves to $\frac{\partial{\hat{H}(\mathbf{s})}}{\partial{s_{i}x}}$ and omit $\frac{\partial{\hat{H}(\mathbf{s})}}{\partials_{i}y}$ to save space):
\begin{equation}
{\label{newweightedparx}}%
\begin{split}
\frac{\partial{\hat{H}(\mathbf{s})}}{\partial{s_{ix}}}=  &  \int_{V(s_{i}%
)}\alpha_{1}(x,\mathbf{s})w_{1}(x,\mathbf{s})\frac{(x-s_{i})_{x}}{d_{i}(x)}dx\\
+  &  \int_{V(s_{i})}\beta_{1}(x,\mathbf{s})\frac{(x-s_{i})_{x}}{d_{i}(x)}dx\\
+  &  \sum_{j\in\Gamma_{i}}sgn(n_{jx})\frac{sin\theta_{ij}}{D_{ij}}\int
_{0}^{z_{ij}}\alpha_{2}(x,\mathbf{s})w_{2}(x,\mathbf{s})rdr\\
+  &  \sum_{j\in\Gamma_{i}}sgn(n_{jx})\frac{sin\theta_{ij}}{D_{ij}}\int
_{0}^{z_{ij}}\beta_{2}(x,\mathbf{s})rdr
\end{split}
\end{equation}
$\frac{\partial{\hat{H}(\mathbf{s})}}{\partial{s_{iy}}}$ can be obtained in a similar way.
Obviously, the boosting process (\ref{modifiedlineari})-(\ref{modifiedlinearb}%
) actually changes the objective function $H(\mathbf{s})$. Thus, when a new
equilibrium is reached in the boosted derivative phase of system operation, it
is necessary to revert to the original objective function by setting
$\alpha_{1}(x,\mathbf{s})=\alpha_{2}(x,\mathbf{s})=1$ and $\beta
_{1}(x,\mathbf{s})=$ $\beta_{2}(x,\mathbf{s})=0$.

We summarize the boosting process as follows. Initially, node $i$ uses
(\ref{weightedparx})-(\ref{weightedpary}) until an equilibrium ${\mathbf{s}%
}^{1}$ is reached at time $\tau^{1}$ and nodes communicate their positions to
each other.

\begin{enumerate}
\item At $t=\tau^{1}$, evaluate $H{(\mathbf{s}(\tau^{1}))}$ and set
${\mathbf{s}}^{\ast}={\mathbf{s}}^{1}$ and $H^{\ast}{=}H{(\mathbf{s}(\tau^{1}%
))}$. Then, apply boosting functions (\ref{modifiedlineari}%
)-(\ref{modifiedlinearb}), evaluate (\ref{newweightedparx}), and iterate on
the controllable node position using (\ref{mov}). Set $BIt=0$. $BIt$ is short for 
the Boosted iteration, which is a counter for iteration needed for a new local optima.

\item Wait until $\frac{\partial{\hat{H}(\mathbf{s})}}{\partial{s_{ix}}}%
=\frac{\partial{\hat{H}(\mathbf{s})}}{\partial{s_{iy}}}=0$ at time $\hat{\tau
}^{1}>\tau^{1}$.

\item At $t=\hat{\tau}^{1}$, set $\alpha_{1}(x,\mathbf{s})=\alpha
_{2}(x,\mathbf{s})=1$ and $\beta_{1}(x,\mathbf{s})=$ $\beta_{2}(x,\mathbf{s}%
)=0$ and revert to $\frac{\partial{H_{i}(\mathbf{s})}}{\partial{s_{i}}}$.

\item Wait until $\frac{\partial{H(\mathbf{s})}}{\partial{s_{ix}}}%
=\frac{\partial H{(\mathbf{s})}}{\partial{s_{iy}}}=0$ at time $\tau^{2}%
>\hat{\tau}^{1}$ and evaluate $H{(\mathbf{s}(\tau^{2}))}$, get $BIt$. If
$H{(\mathbf{s}(\tau^{2}))>}H^{\ast}$, then set ${\mathbf{s}}^{\ast
}={\mathbf{s}(\tau^{2})}$ and $H^{\ast}{=}H{(\mathbf{s}(\tau^{2}))}$.
Otherwise, ${\mathbf{s}}^{\ast},$ $H^{\ast}$ remain unchanged (if nodes are
mobile and have already been moved to ${\mathbf{s}(\tau^{2})}$, then return
them to ${\mathbf{s}}^{\ast}$).

\item Either STOP, or repeat the process from the current ${\mathbf{s}}^{\ast
}$ with a new boosting function to further explore the mission space for
better equilibrium points.
\end{enumerate}

Note that if ${\mathbf{s}}^{1}$ is a global optimum, then the boosting process
simply perturbs node locations until Step 4  returns them to
${\mathbf{s}}^{1}$. The process will stop if no solution is better than $\mathbf{s}^{1}$
 after trying finite boosting functions.
 It is also possible (due to symmetry) that there are
multiple global optima, in which case $H{(\mathbf{s}(\tau^{2}))=}%
H{(\mathbf{s}(\tau^{1}))}$ and the new equilibrium point is equivalent to the
original one. 

The process above assumes that all nodes wait until they have all reached an
equilibrium point ${\mathbf{s}}^{1}$ before each initiates its boosting
process. However, this may also be done in a distributed function through a
self-boosting process: node $i$ may apply (\ref{modifiedlineari}%
)-(\ref{modifiedlinearb}) as soon as it observes $\frac{\partial
{H(\mathbf{s})}}{\partial{s_{ix}}}=\frac{\partial H{(\mathbf{s})}}%
{\partial{s_{iy}}}=0$.

\subsection{Boosting Function Selection}

The selection of boosting functions generally depends on the mission space
topology. For instance, it is clear that if there are no obstacles, then
$\alpha_{2}(x,\mathbf{s})=1$, $\beta_{2}(x,\mathbf{s})=0$, since only the
first integrals in (\ref{weightedparx})-(\ref{weightedpary}) are present. In
what follows, we present three families of boosting functions that we have
investigated to date; each has different properties and has provided promising results.

Before proceeding, we make a few observations which guide the selection of
boosting functions. First, we exclude cases such that $\alpha_{1}%
(x,s_{i})=\alpha_{2}(x,s_{i})=C$ independent of $x$, and $\beta_{1}(x,s_{i})=$
$\beta_{2}(x,s_{i})=0$. In such cases, the boosting effect is null, since it
implies that $\frac{\partial{\hat{H}(\mathbf{s})}}{\partial{s_{i}}}%
=C\frac{\partial{H(\mathbf{s})}}{\partial{s_{i}}}$, which has no effect on
$\frac{\partial{H(\mathbf{s})}}{\partial{s_{i}}}=0$. Second, we observe that
if $\left\vert \beta_{1}(x,s_{i})\right\vert >>\alpha_{1}(x,s_{i}%
)w_{1}(x,\mathbf{s})$, then the first integral in (\ref{newweightedparx}) is
dominated by the second one, and the net effect is that nodes tend to be
attracted to a single point (their center of mass) instead of exploring the
mission space. The third observation is more subtle. The first term of
(\ref{weightedparx}) contains information on points of the visible set
$V(s_{i})$, which is generally more valuable (i.e., more points in $V(s_{i})$)
than the information in the second term related to the boundary points in
$\Gamma_{i}$ (except, possibly, for unusual obstacle configurations). Thus, a
boosting function should ensure that the first integral in (\ref{weightedparx}%
) dominates the second when $\frac{\partial{H_{i}(\mathbf{s})}}{\partial
{s_{i}x}}\neq0$. In order to avoid such issues, in the sequel we limit
ourselves to boosting $w_{1}(x,\mathbf{s})$ only and, therefore, we set $\alpha
_{2}(x,s_{i})=1$, $\beta_{2}(x,s_{i})=0$.

\subsubsection{$P$-Boosting function}

In this function, we keep $\beta_{1}(x,\mathbf{s})=0$ and only concentrate on
$\alpha_{1}(x,\mathbf{s})$ which we set:
\begin{equation}
\alpha_{1}(x,\mathbf{s})=k{P(x,\mathbf{s})}^{-\gamma}%
\end{equation}
where $P(x,\mathbf{s})$ is the joint detection probability defined in
(\ref{jointP}), $\gamma$ is a positive integer parameter and $k$ is a gain
parameter. Thus, the boosted derivative associated with this $P$-boosting
function is
\begin{equation}
{\label{power boost}}
\begin{split}
\frac{\partial{\hat{H}(\mathbf{s})}}{\partial{s_{ix}}}=  &  \int_{V(s_{i}%
)}k{P(x,\mathbf{s})}^{-\gamma}w_{1}(x,\mathbf{s})\frac{(x-s_{i})_{x}}{d_{i}(x)}dx\\
+  &  \sum_{j\in\Gamma_{i}}sgn(n_{jx})\frac{sin\theta_{ij}}{D_{ij}}\int
_{0}^{z_{ij}}w_{2}(x,\mathbf{s})rdr\\
\end{split}
\end{equation}
The motivation for this function is similar to a method used in
\cite{Minyi2011} to assign higher weights for low-coverage interior points in
$V(s_{i})$, in order for nodes to explore such low coverage areas. This is
consistent with the following properties of this boosting function:
$({P(x,\mathbf{s})})^{-\gamma}\rightarrow\infty$ as $P(x,\mathbf{s}%
)\rightarrow0$, and $({P(x,\mathbf{s})})^{-\gamma}\rightarrow1$ as
$P(x,\mathbf{s})\rightarrow1$.

\subsubsection{Neighbor-Boosting function}

We set $\alpha_{1}(x,\mathbf{s})=1$ and focus on $\beta_{1}(x,\mathbf{s})$.
Every node applies a repelling force on each of its neighbors with the effect
being monotonically decreasing with their relative distance. We define:
\begin{equation}
\beta_{1}(x,\mathbf{s})=\sum_{j\in B_{i}}\delta(x-s_{j})\frac{k_{j}}{\Vert
s_{i}-x\Vert^{\gamma}}%
\end{equation}
where $k_{j}\geq0$ is a gain parameter for $j$, $\gamma$ is a positive integer
parameter, and $\delta(x-s_{j})$ is the delta function. The boosted derivative
associated with this neighbor-boosting function is
\begin{equation}
{\label{neighbor boost}}%
\begin{split}
\frac{\partial{\hat{H}(\mathbf{s})}}{\partial{s_{ix}}}= &  \int_{V(s_{i}%
)}w_{1}(x,\mathbf{s})\frac{(x-s_{i})_{x}}{d_{i}(x)}dx\\
+ &  \sum_{j\in\Gamma_{i}}sgn(n_{jx})\frac{sin\theta_{ij}}{D_{ij}}\int
_{0}^{z_{ij}}w_{2}(x,\mathbf{s})rdr\\
+ &  \sum_{j\in B_{i}}\frac{k_{j}}{\Vert s_{j}-s_{i}\Vert^{\gamma+1}}%
(s_{j}-s_{i})_{x}%
\end{split}
\end{equation}
Note that $k_{j}$ may vary over different neighbors $j$. For instance, if some
neighboring node $j$ is such that $j\notin V(s_{i})$, then we may set
$k_{j}=0$.
\subsubsection{$\Phi$-boosting function}

This function aims at varying $\alpha_{1}(x,\mathbf{s})$ by means of $\Phi
_{i}(x)$ defined in (\ref{phi}), which is the probability that point $x$ is
not detected by neighboring nodes of $i$. $\beta_{1}(x,\mathbf{s})=0$ as well.
Large $\Phi_{i}(x)$ values imply a
lower coverage by neighbors, therefore higher weights are set. In particular,
we define
\begin{equation}
\alpha_{1}(x,\mathbf{s})=k\Phi_{i}(x)^{\gamma}%
\end{equation}
where $k$ is a gain parameter and $\gamma$ is a positive integer parameter.
The boosted derivative here is
\begin{equation}
{\label{phi boost}}%
\begin{split}
\frac{\partial{\hat{H}(\mathbf{s})}}{\partial{s_{ix}}}=  &  \int_{V(s_{i}%
)}k\Phi_{i}(x)^{\gamma}w_{1}(x,\mathbf{s})\frac{(x-s_{i})_{x}}{d_{i}(x)}dx\\
+  &  \sum_{j\in\Gamma_{i}}sgn(n_{jx})\frac{sin\theta_{ij}}{D_{ij}}\int
_{0}^{z_{ij}}w_{2}(x,\mathbf{s})rdr\\
&
\end{split}
\end{equation}
Observe that $\Phi_{i}(x)=0$ means that $x$ is well-covered by neighbors of
$i$, therefore, sensor node $i$ has no incentive to move closer to this point.
On the other hand, $\Phi_{i}(x)=1$ means that no neighbor covers $x$, so the
boosted weight is the value of the gain $k$.

To compare the performance of the boosting function method to that of a random perturbation method, 
we propose a random perturbation method applied in step 1 to get (\ref{newweightedparx}) in the boosting process.
Let $\xi_x, \xi_y$ be independent random variables. 
The perturbed derivatives $\frac{\partial{\hat{H}(\mathbf{s})}}{\partial{s_{ix}}}$, $\frac{\partial{\hat{H}(\mathbf{s})}}{\partial{s_{iy}}}$ will be
\begin{equation}{\label{randomperturbx}}
\frac{\partial{\hat{H}(\mathbf{s})}}{\partial{s_{ix}}} = \frac{\partial{H(\mathbf{s})}}{\partial{s_{ix}}} + \xi_x 
\end{equation}
\begin{equation}{\label{randomperturby}}
\frac{\partial{\hat{H}(\mathbf{s})}}{\partial{s_{iy}}} = \frac{\partial{H(\mathbf{s})}}{\partial{s_{iy}}} + \xi_y 
\end{equation}
Note that $\xi_x$ and $\xi_y$ are independently updated for each node in each iteration.
Then, this random perturbation method can be performed in a distributed way. 
\section{Simulation Results}

In this section, we provide simulation examples illustrating how the objective
function value in (\ref{obj}) is improved by using the boosting function
process and how the parameter values in the  boosting functions we have
considered can further affect performance. Moreover, 
we show how the boosting method is superior to the random perturbation approach in terms of
the number of iterations to a new local optimum.

\begin{figure}[h]%
\begin{tabular}
[c]{cc}%
\begin{minipage}[t]{1.5in} \includegraphics[
	height=1.25in,
	width=1.45in]
	{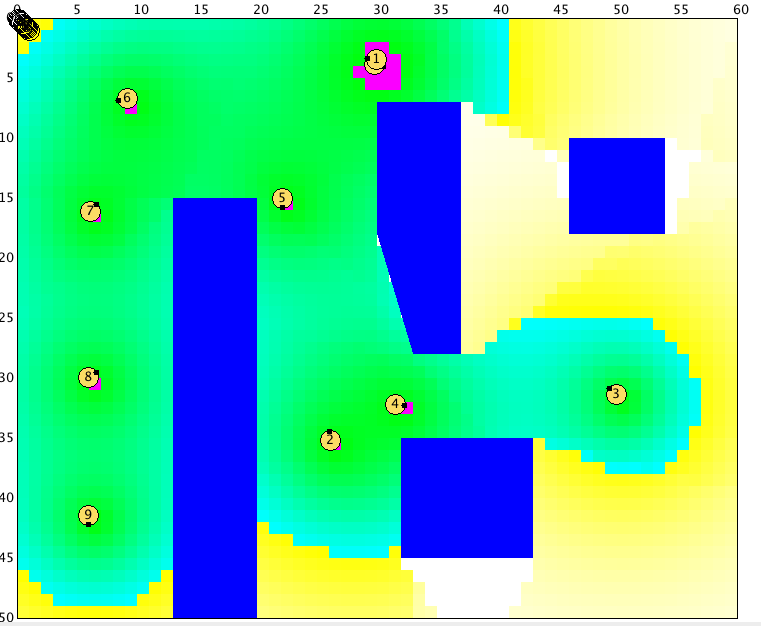} \subcaption{General obstacle with \\$H(\textbf{s}_0^*)=1368.3$} \label{figGeneralInitial} \end{minipage}\begin{minipage}[t]{1.5in} \includegraphics[
	height=1.25in,
	width=1.45in]
	{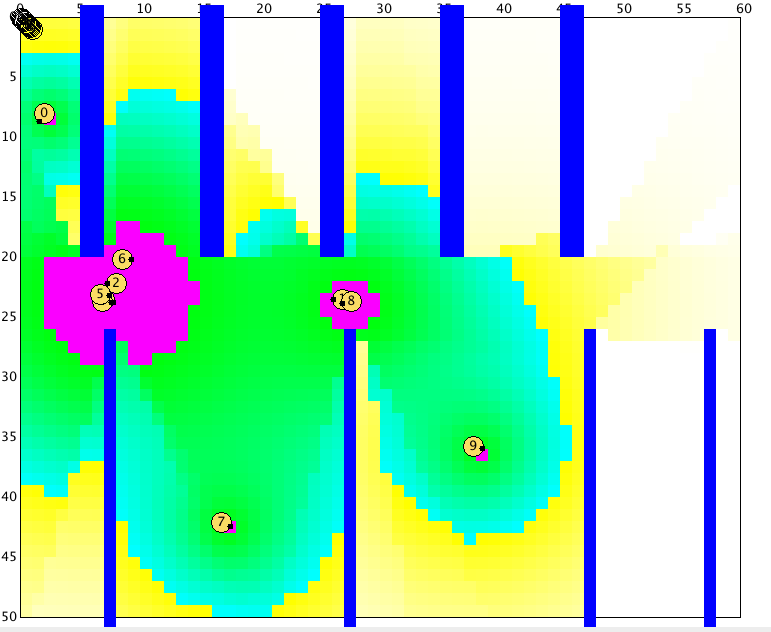} \subcaption{Room obstacle with \\ $H(\textbf{s}_0^*)=1183.5$} \label{figRoomInitial} \end{minipage} &
\\
\begin{minipage}[t]{1.5in} \includegraphics[
	height=1.25in,
	width=1.45in]
	{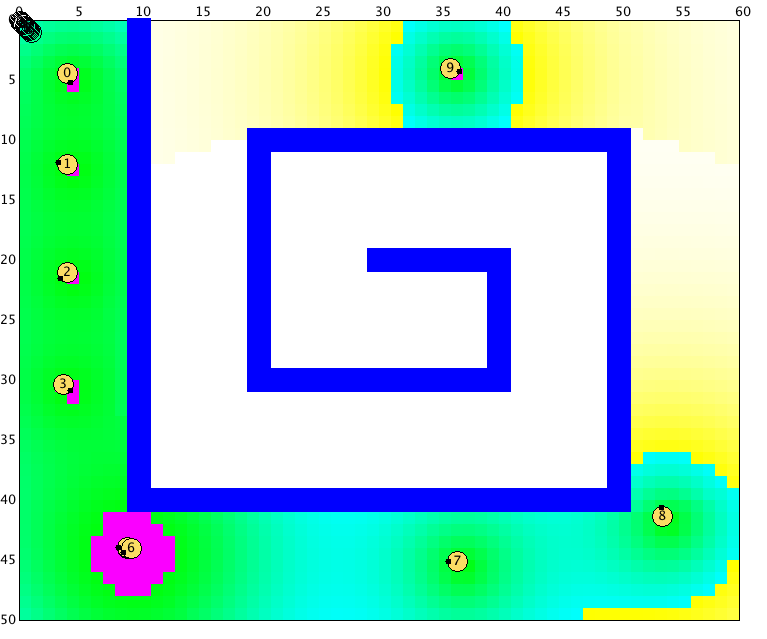} \subcaption{Maze obstacle with \\ $H(\textbf{s}_0^*)=904.0$} \label{figMazeInitial} \end{minipage}\begin{minipage}[t]{1.5in} \includegraphics[
	height=1.25in,
	width=1.45in]
	{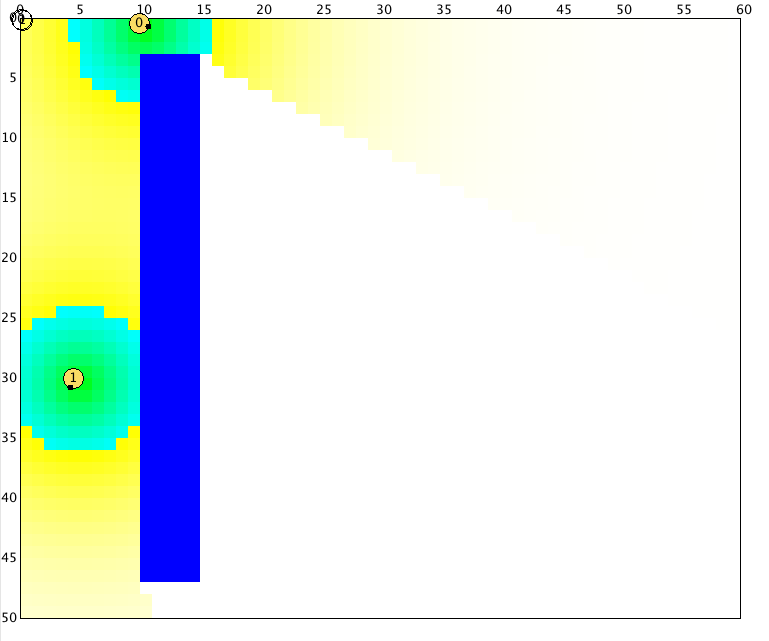} \subcaption{Narrow obstacle with \\$H(\textbf{s}_0^*)=246.5$} \label{figNarrowInitial}\end{minipage} &
\end{tabular}
\caption{Initial local optima in all obstacle configurations}
\label{figInitial}
\end{figure}

%Figures \ref{figGeneralInitial}-\ref{figNarrowInitial} present four mission
Figure. \ref{figInitial} presents four mission
spaces with different obstacle configurations (obstacles shown as blue
polygons), which we refer to as \textquotedblleft\textit{General}
Obstacle\textquotedblright, \textquotedblleft\textit{Room}
Obstacle\textquotedblright, \textquotedblleft\textit{Maze}
Obstacle\textquotedblright\ and \textquotedblleft\textit{Narrow}
Obstacle\textquotedblright\,, respectively. The event density functions are uniform in 
all cases, i.e., $R(x)=1$. In the first three cases, there are
10 nodes shown as numbered circles while in the \textit{Narrow} Obstacle case, there are only 2 nodes.
The mission space is colored from dark to lighter as the joint detection
probability decreases (the joint detection
probability is $\geq0.97$ for purple areas, $\geq0.50$ for green areas, and
near zero for white areas). Nodes start from the upper left corner and reach
equilibrium configurations obtained by the gradient-based algorithm in \cite{Minyi2011}.
The objective function values at the equilibria are shown in
the captions of Figs. \ref{figGeneralInitial}-\ref{figNarrowInitial}. 
It is easy to see that these deployments are sub-optimal due to
the obvious imbalanced coverage. For instance, in Fig. \ref{figRoomInitial},
the upper and lower rightmost \textquotedblleft rooms\textquotedblright\ are
poorly covered while there are 4 nodes clustered together near the first
obstacle on the left side. We expect that boosting functions can guide nodes
towards exploration of poorly covered areas in the mission space, thus leading
to a more balanced, possibly globally optimal, equilibrium. 

First, we discuss how we select parameters for the boosting functions.
For the neighbor-boosting function, we select the gain
parameters $k_{j}$ in two different ways: $(i)$ the same
for all neighboring nodes in a line of sight of $s_{i}$, otherwise,
$k_{j}=0$:
\begin{equation}
k_{j}=%
\begin{cases}
k & \text{if}\quad s_{j}\in V(s_{i}),\text{ }j\in B_{i}\\
0 & \text{otherwise}%
\end{cases}
\end{equation}
and $(ii)$, $k_{j}=0$ for all neighboring nodes except for the closest
neighbor of $s_{i}$:
\begin{equation}
k_{j}=%
\begin{cases}
k & j=\text{arg }min_{k\in B_{i}}\Vert s_{i}-s_{k}\Vert\\
0 & \text{otherwise}%
\end{cases}
\label{kj-2}%
\end{equation}
We define $H(\mathbf{s}^{\ast})_{1}$ and $H(\mathbf{s}^{\ast})_{2}$ to
correspond to the objective function values after the boosting process for
each of these two choices and have found through extensive experimentation
(shown in Table. \ref{table2}) that $H(\mathbf{s}^{\ast})_{2}>H(\mathbf{s}^{\ast})_{1}$ for
almost cases considered. In the following discussion, the second definition of $k_{j}$ is used.

We also study the effect of the parameters $\gamma$ and $k$
and have found the $\gamma$, $k$ that yield the best
results for all boosting functions (shown in the captions). Table \ref{table2} lists results from some of our experiments. For instance, in the room case, the neighbor-boosting function with $\gamma=1$ and $k=300$ yields the largest objective value $H(\mathbf{s}^{\ast})_{2}$. 

\begin{figure}[h]%
\begin{tabular}
[c]{cc}%
\begin{minipage}[t]{1.5in} \includegraphics[
	height=1.25in,
	width=1.45in]
	{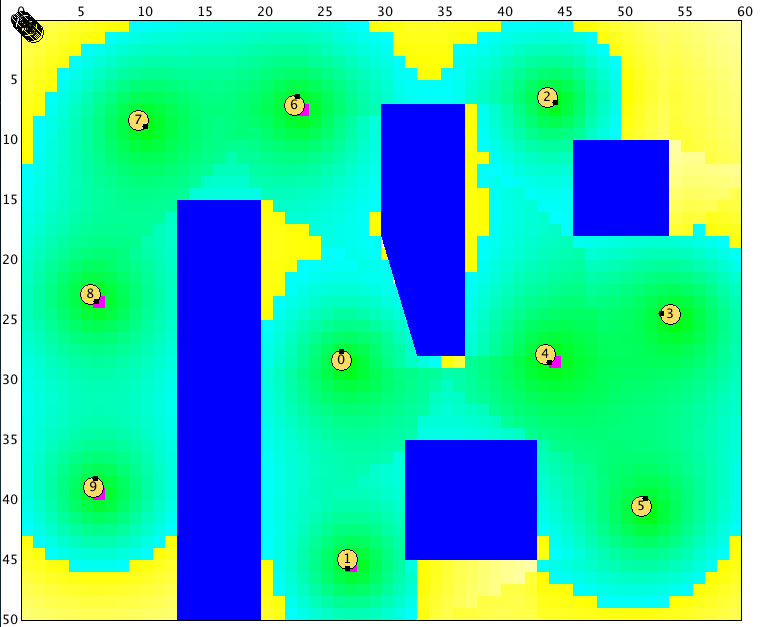} \subcaption{$P$-boost, $\gamma=4$, $k = 100$,\\BIt=161; $H(\textbf{s}^*)=1533.6$} \label{figpower_general} \end{minipage}\begin{minipage}[t]{1.5in} \includegraphics[
	height=1.25in,
	width=1.45in]
	{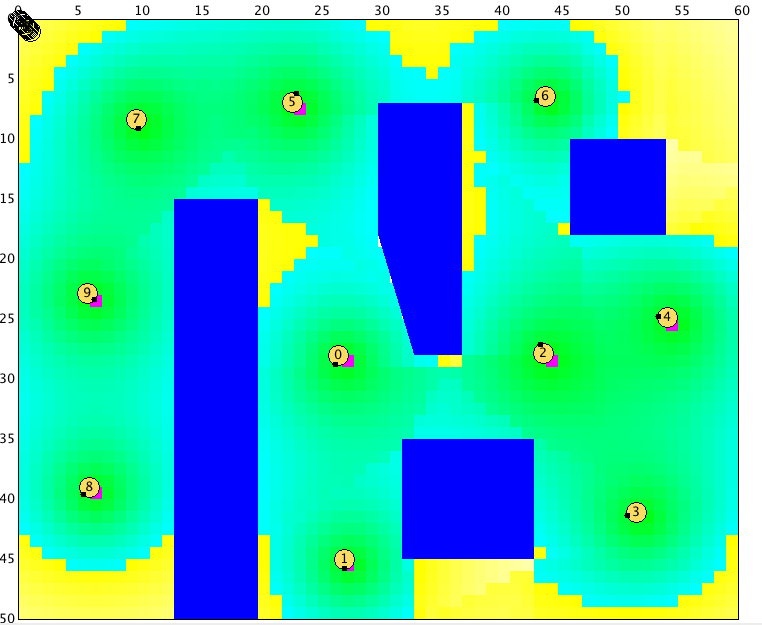} \subcaption{Neighbor-boost, $\gamma=2$,\\ $k = 500$, BIt=390; $H(\textbf{s}^*)=1533.3$} \label{figneighbor_general} \end{minipage} &
\\
\begin{minipage}[t]{1.5in} \includegraphics[
	height=1.25in,
	width=1.45in]
	{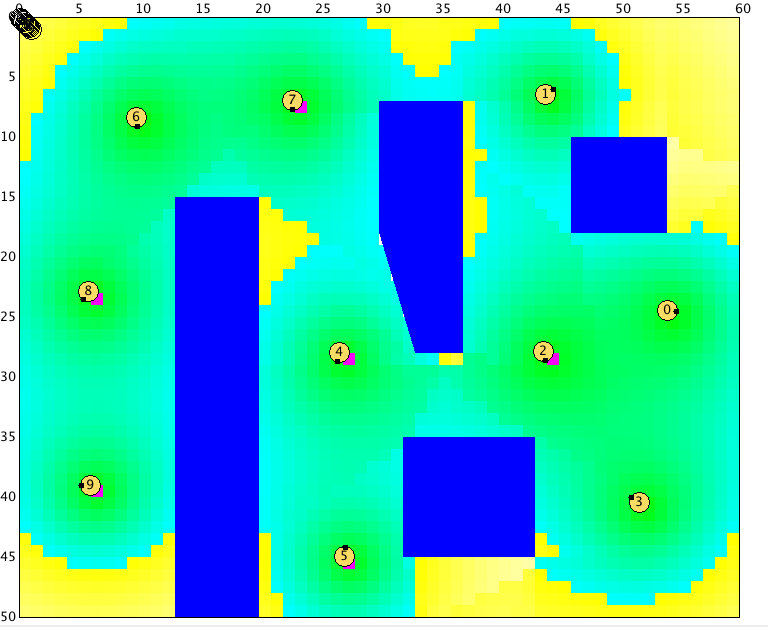} \subcaption{$\Phi$-boost, $\gamma=2$, $k=1000$, \\ BIt=160; $H(\textbf{s}^*)=1533.4$} \label{figphi_general} \end{minipage}\begin{minipage}[t]{1.5in} \includegraphics[
	height=1.25in,
	width=1.45in]
	{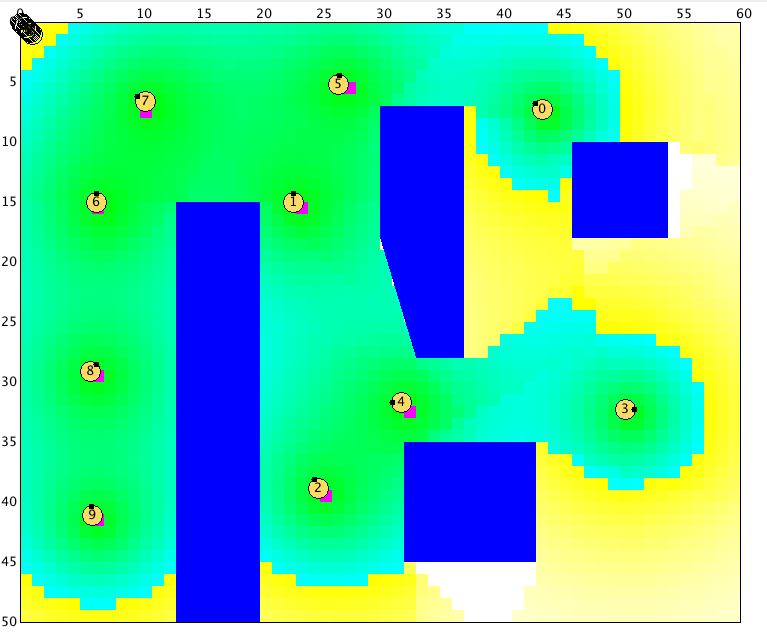} \subcaption{Random Perturbation, \\ BIt=653; $H(\textbf{s}^*)=1443.9$} \label{figrandom_general} \end{minipage} &
\end{tabular}
\caption{General Obstacle Configuration}
\label{figGeneral}
\end{figure}

Then, we show the results for all configurations.
%Figures \ref{figpower_general}-\ref{figrandom_general} illustrate the effects of
Figure. \ref{figGeneral} illustrates the effects of
different methods used in the general obstacle configuration.
The $P$-boosting and the $\Phi$-boosting functions attain the best local optima (objective values are increased by 12\%) in the smallest number of iterations. Figure. \ref{figrandom_general} shows a snapshot of a typical result using the random perturbation approach in (\ref{randomperturbx})-(\ref{randomperturby}). It needs about four times as many iterations as the $\Phi$-boosting function, yet converges to a worse local optimum.

\begin{figure}[h]%
\begin{tabular}
[c]{cc}%
\begin{minipage}[t]{1.5in} \includegraphics[
	height=1.25in,
	width=1.45in]
	{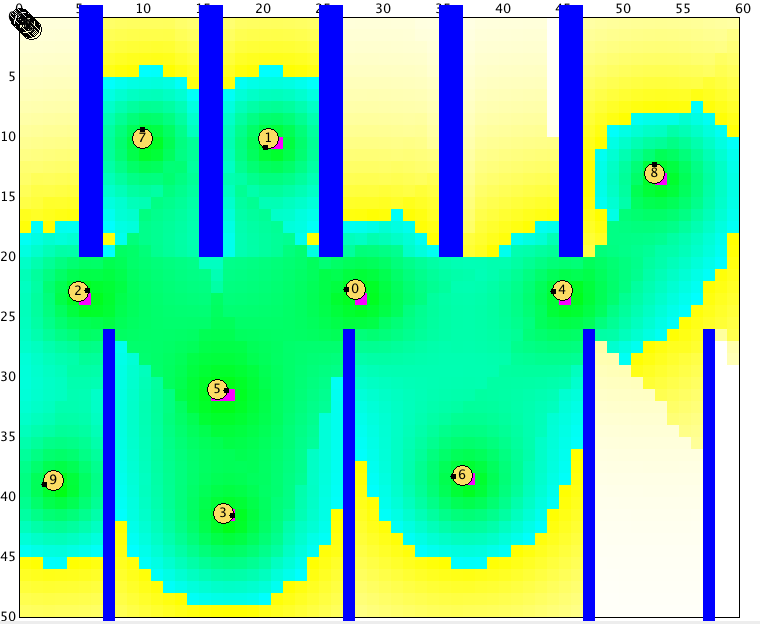} \subcaption{$P$-boost, $\gamma=4$, $k = 100$; \\ BIt=221, $H(\textbf{s}^*)=1419.5$} \label{figpower_room} \end{minipage} \begin{minipage}[t]{1.5in} \includegraphics[
	height=1.25in,
	width=1.45in]
	{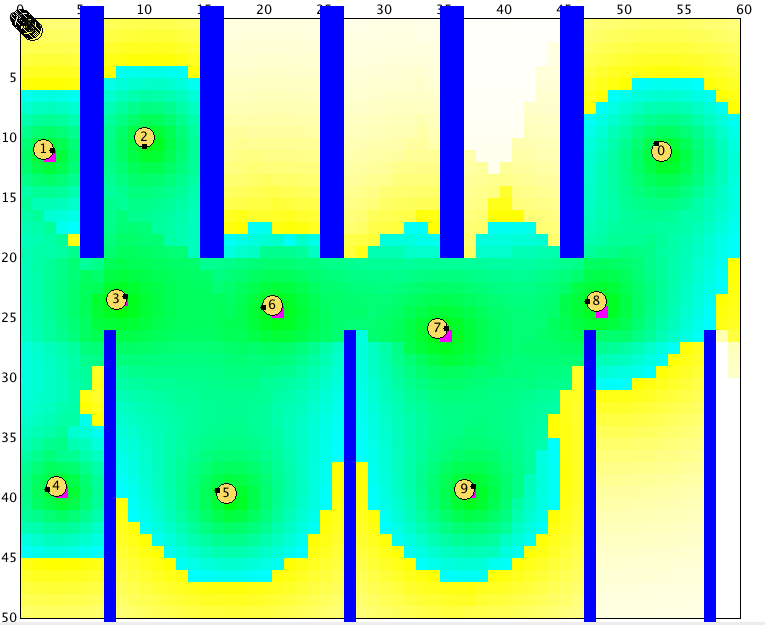} \subcaption{Neighbor-boost, $\gamma=1$, \\ $k = 300$; BIt=364, $H(\textbf{s}^*)=1417.1$} \label{figneighbor_room} \end{minipage} &
\\
\begin{minipage}[t]{1.5in} \includegraphics[
	height=1.25in,
	width=1.45in]
	{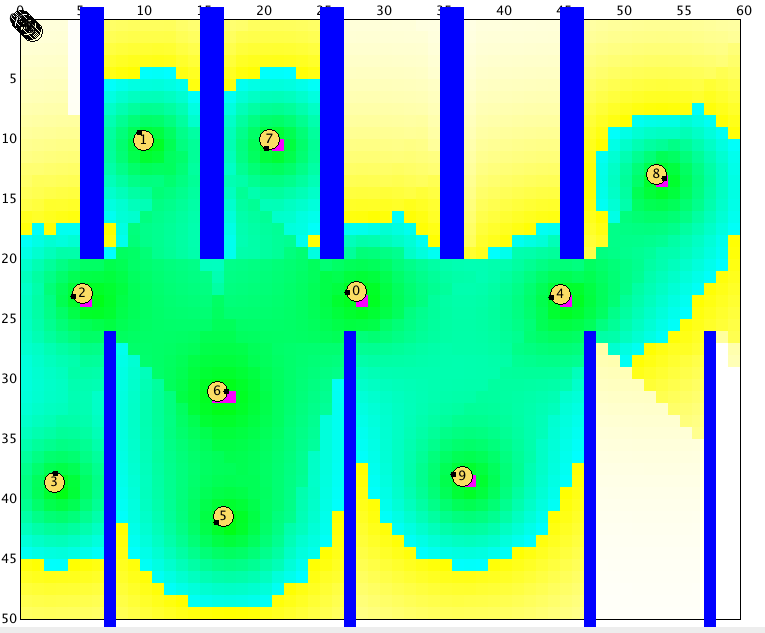} \subcaption{$\Phi$-boost, $\gamma=1$, $k=1000$;\\ BIt=208, $H(\textbf{s}^*)=1419.1$} \label{figphi_room} \end{minipage} \begin{minipage}[t]{1.5in} \includegraphics[
	height=1.25in,
	width=1.45in]
	{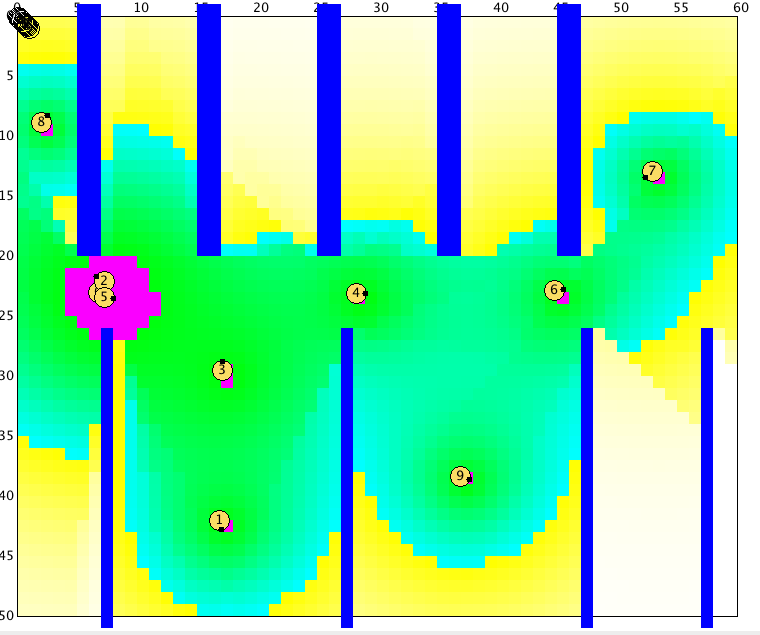} \subcaption{Random Perturbation \\ BIt=850; $H(\textbf{s}^*)=1377.3$} \label{figrandom_room} \end{minipage} &
\end{tabular}
\caption{Room Obstacle Configuration}
\label{figRoom}
\end{figure}

%Next, we consider the "Room" obstacle case. Fig. \ref{figpower_room}-\ref{figrandom_room} show similar
Next, we consider the "Room" obstacle case in Fig. \ref{figRoom}.
Comparing Fig. \ref{figRoomInitial} with Fig. \ref{figRoom}, the clustered nodes in Fig. \ref{figRoomInitial} have spread apart and the objective value has increased. The $P$-boosting and the $\Phi$-boosting converge to better local optima (about 20\% increase in the objective function value over the original one) than those resulting from the neighbor-boosting function. The random perturbation gets stuck at a worse equilibrium after more iterations than any boosting function. 

\begin{figure}[h]%
\begin{tabular}
[c]{cc}%
\begin{minipage}[t]{1.5in} \includegraphics[
	height=1.25in,
	width=1.45in]
	{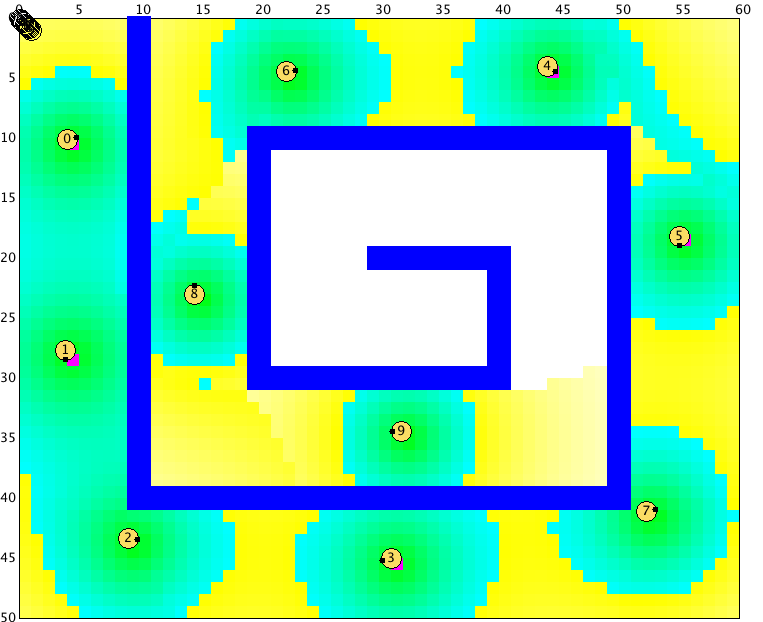} \subcaption{$P$-boost, $\gamma=4$, $k = 100$; \\ BIt=517, $H(\textbf{s}^*)=1180.5$} \label{figpower_maze} \end{minipage}\begin{minipage}[t]{1.5in} \includegraphics[
	height=1.25in,
	width=1.45in]
	{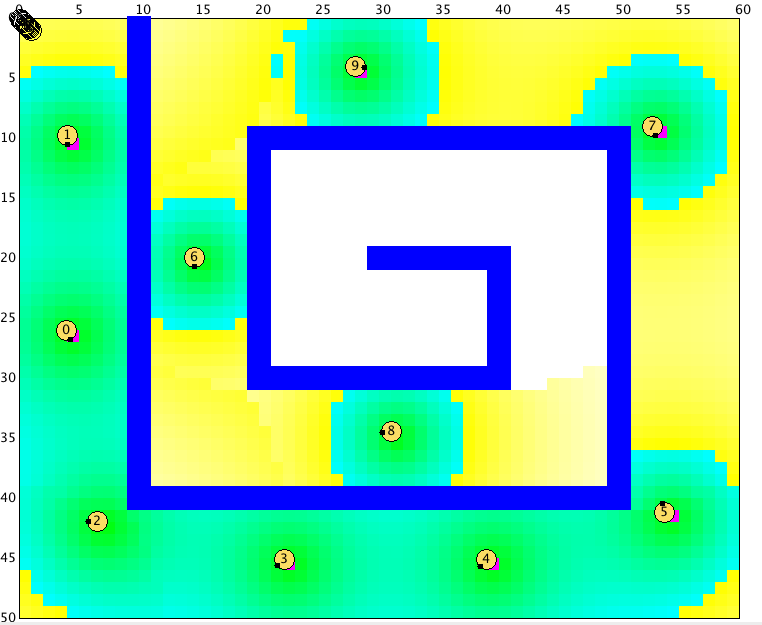} \subcaption{Neighbor-boost,$\gamma=2$, \\ $k = 1000$; BIt=600, $H(\textbf{s}^*)=1168.6$} \label{figneighbor_maze} \end{minipage} &
\\
\begin{minipage}[t]{1.5in} \includegraphics[
	height=1.25in,
	width=1.45in]
	{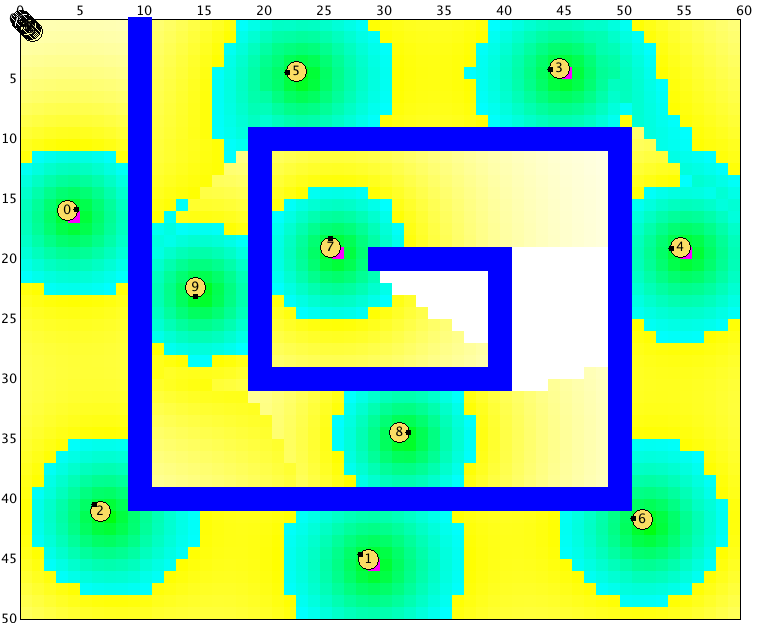} \subcaption{$\Phi$-boost, $\gamma=2$, $k=100$;\\ BIt=503, $H(\textbf{s}^*)=1236.1$} \label{figphi_maze} \end{minipage} \begin{minipage}[t]{1.5in} \includegraphics[
	height=1.25in,
	width=1.45in]
	{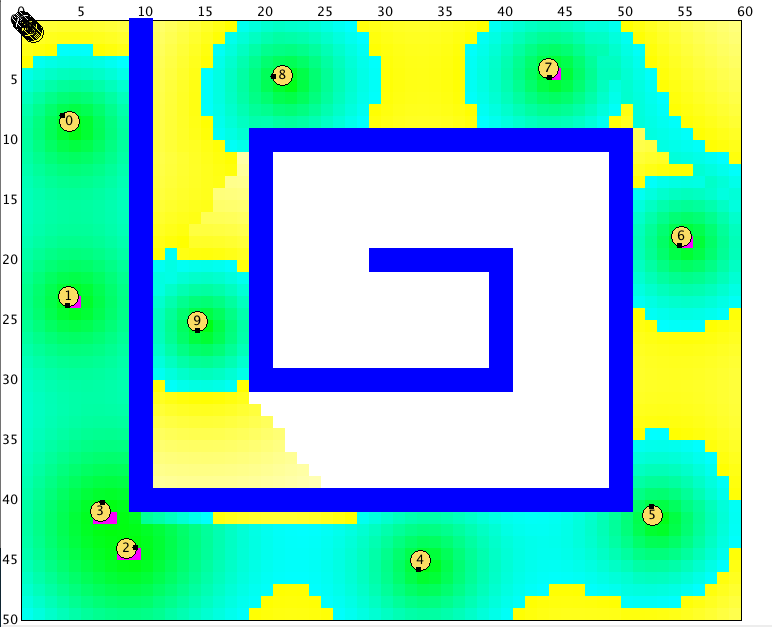} \subcaption{Random Perturbation \\ BIt=3439; $H(\textbf{s}^*)=1132.3$} \label{figrandom_maze} \end{minipage} &
\end{tabular}
\caption{Maze Obstacle Configuration}
\label{figMaze}
\end{figure}

\begin{figure}[h]%
\begin{tabular}
[c]{cc}%
\begin{minipage}[t]{1.5in} \includegraphics[
	height=1.25in,
	width=1.45in]
	{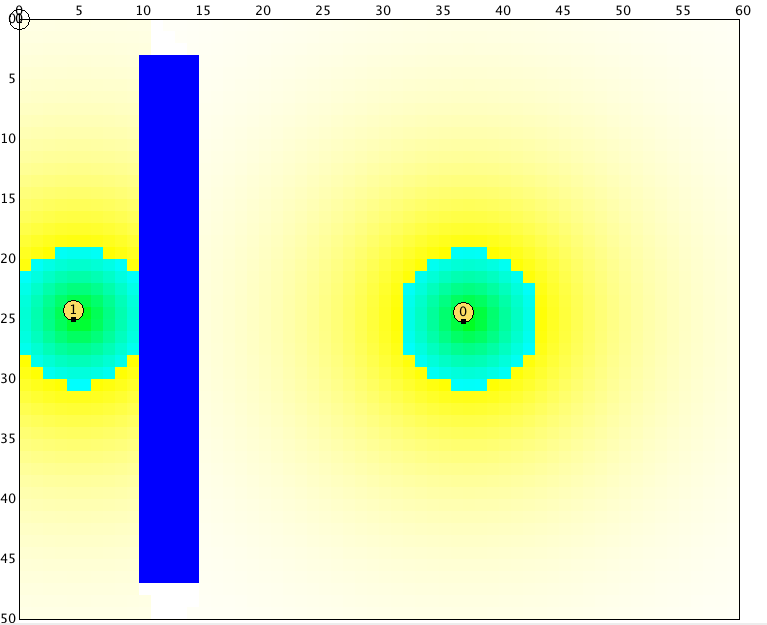} \subcaption{$P$-boost, $\gamma=4$, $k = 100$;\\ BIt=103, $H(\textbf{s}^*)=502.5$} \label{figpower_narrow} \end{minipage} \begin{minipage}[t]{1.5in} \includegraphics[
	height=1.25in,
	width=1.45in]
	{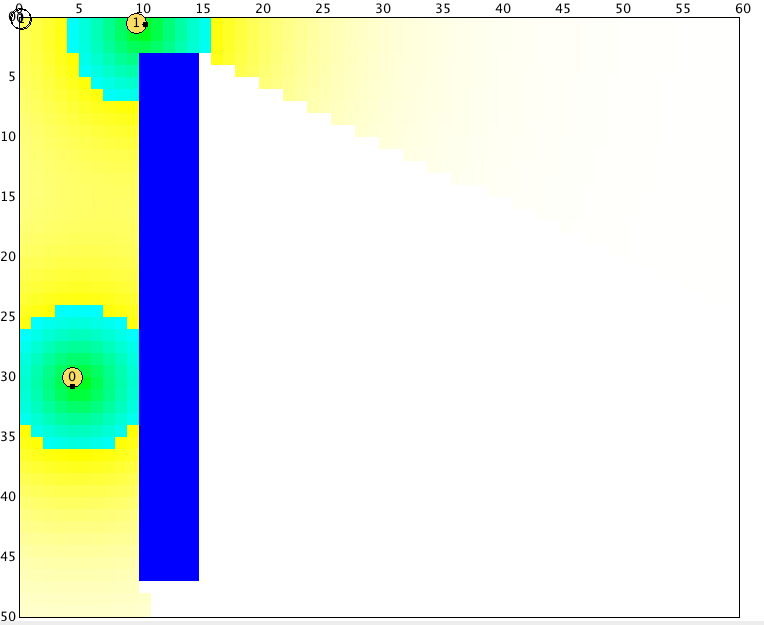} \subcaption{Neighbor-boost,  $\gamma=1$, \\$k = 300$; BIt=212; $H(\textbf{s}^*)=246.5$} \label{figneighbor_narrow} \end{minipage}  &
\\
\begin{minipage}[t]{1.5in} \includegraphics[
	height=1.25in,
	width=1.45in]
	{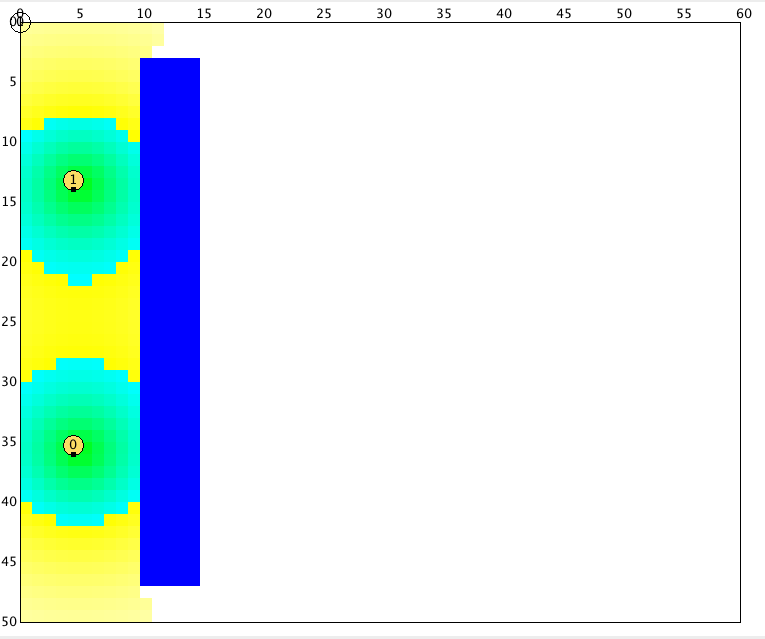} \subcaption{$\Phi$-boost, $\gamma=2$, $k=1000$; \\ BIt=90; $H(\textbf{s}^*)=253.3$} \label{figphi_narrow} \end{minipage} 
	\begin{minipage}[t]{1.5in} \includegraphics[
	height=1.25in,
	width=1.45in]
	{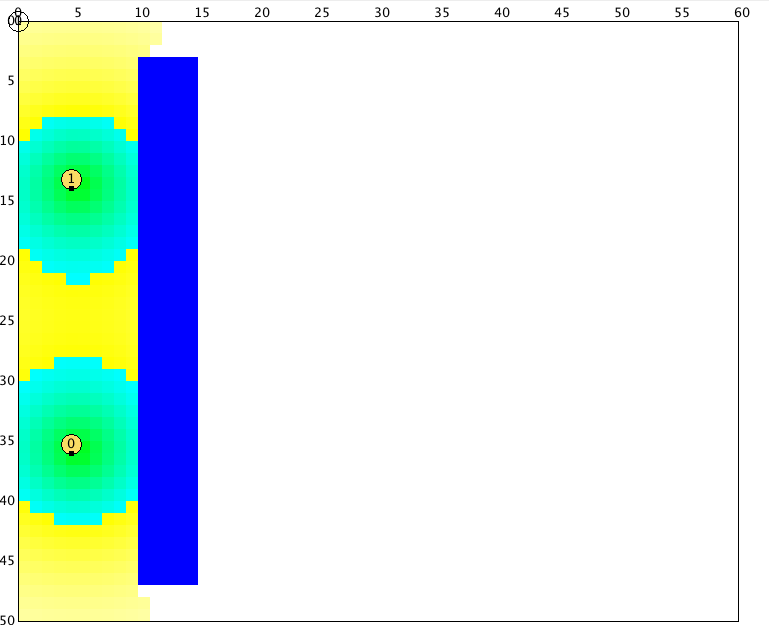} \subcaption{Random Perturbation \\ BIt=1153; $H(\textbf{s}^*)=253.3$} \label{figrandom_narrow} \end{minipage} &
\end{tabular}
\caption{Narrow Obstacle Configuration}
\label{figNarrow}
\end{figure}

%Figures \ref{figpower_maze}-\ref{figrandom_maze} display the results of boosting functions applied to the
Figure. \ref{figMaze} displays the results of boosting functions applied to the maze configuration. The $\Phi$-boosting function attains a local optimum with the highest objective function value (approximately a 44\% increase in the objective function value over the original one) among all methods while the random-boosting does the worst. 
%Figures \ref{figpower_narrow}-\ref{figrandom_narrow} 
Figure. \ref{figNarrow}
shows results for the narrow obstacle configuration where the $P$-boosting function 
works the best and the objective function value is increased by 105\%, from 245.3 to 502.5. Note that the neighbor-boosting function fails to escape the local optimum. This is because the repelling forces between the two nodes have no components to drive sensor nodes over the obstacle. 
Although the random perturbation method converges to similar results as the $\Phi$-boosting function, it requires many more iterations.

In summary,  we conclude that the boosting function approach, while still not
guaranteeing global optimality, provides substantial improvements in the
objective function value, varying from 12\% to 105\%. In addition,  
the boosting function approach converges to an equilibrium faster and usually with a higher objective function value than the random perturbation method.

%\begin{table}[th]
%\renewcommand{\arraystretch}{1.5}
%\caption{Summary of boosting function effects}%
%\centering%
%\begin{tabular}
%[c]{|c|c|c|c|c|}\hline
%Space & Best boosting & $H(\mathbf{s}^{\ast})$ & $H(\mathbf{s}_{0}^{\ast})$ &
%$\frac{H(\mathbf{s}^{\ast})-H(\mathbf{s}_{0}^{\ast})}{H(\mathbf{s}_{0}^{\ast
%})}$\\\hline
%General & All & 1530.3 & 1368.3 & 12\%\\
%Room & $P$ and $\Phi$ & 1419.1 & 1183.5 & 20\%\\
%Maze & $\Phi$ & 1236.1 & 860.7 & 44\%\\
%Narrow & $P$ and $\Phi$ & 502.5 & 245.5 & 105\%\\\hline
%\end{tabular}
%\label{table4}%
%\end{table}

\begin{table}[ht]
\centering
\begin{tabular}{|c|c|c|c|c|}
\hline
 $\gamma$ & $k$ & Obstacle Type & $H(\textbf{s}^*)_1$ & $H(\textbf{s}^*)_2$\\ 
 
 \hline
 1 & 300 & General & 1513.7 & 1470.0\\
 2 & 300 & General & 1450.6 & 1451.0\\
 2 & 500 & General & 1505.1 & 1533.3\\
 2 & 1000& General & 1446.6 & 1530.7\\
 \hline
 1 & 300 & Room & 1372.9 & 1417.1 \\
 2 & 300 & Room & 1380.8 & 1392.5 \\
 2 & 500 & Room & 1382.8 & 1395.2 \\
 2 & 1000 & Room & 1378 &  1380.8\\
 \hline
 1 & 300 & Maze & 1051.8 & 1110.3\\
 2 & 300 & Maze & 1051.8 & 1133.7\\
 2 & 500 & Maze & 1109.3 & 1110.5 \\
 2 & 1000 & Maze & 1133.9 & 1168.6\\
 \hline
 1 & 300 & Narrow & 245.3 & 245.3\\
 2 & 300 & Narrow & 245.3 & 245.3\\
 2 & 500 & Narrow & 245.3 & 245.3\\
 2 & 1000 & Narrow & 245.3 & 245.3\\
 \hline 
 
\end{tabular}
\caption{The boosted objective function values by neighbor-boosting}
\label{table2}
\end{table}

\section{Conclusions and future work}

\label{sec:concl} We have shown that the objective
function $H(\mathbf{s})$ for the class of optimal coverage control
problems in multi-agent system environments can be decomposed into a local
objective function $H_{i}(\mathbf{s})$ for each node $i$ and a function
independent of node $i$'s controllable position $s_{i}$. This leads to the
definition of boosting functions to systematically (as opposed to randomly)
allow nodes to escape from a local optimum so that the attraction exerted by
some points on \ a node $i$ is \textquotedblleft boosted\textquotedblright\ to
promote exploration of the mission space by $i$ in search of better optima. We
have defined three families of boosting functions, and provided simulation
results illustrating their effects and relative performance. Ongoing research
aims at combining different boosting functions to create a \textquotedblleft
hybrid\textquotedblright\ approach and at studying
self-boosting processes whereby individual nodes autonomously control their
boosting in a distributed manner.

\bibliographystyle{IEEEtran}
\bibliography{research}

\end{document}